# GLOBAL WELL-POSEDNESS FOR THE FOCUSING CUBIC NLS ON THE PRODUCT SPACE $\mathbb{R} \times \mathbb{T}^3$

XUEYING YU, HAITIAN YUE, AND ZEHUA ZHAO

ABSTRACT. In this paper, we prove the global well-posedness for the focusing, cubic nonlinear Schrödinger equation on the product space $\mathbb{R} \times \mathbb{T}^3$ with initial data below the threshold that arises from the the ground state in the Euclidean setting. The defocusing analogue was discussed and proved in Ionescu-Pausader [22] (Comm. Math. Phys. 312 (2012), no. 3, 781-831).

## 1. Introduction

We study the long time behavior of the focusing, cubic nonlinear Schrödinger equation (NLS) on $\mathbb{R} \times \mathbb{T}^3$ with initial data in the energy space $H^1$

$$\begin{cases} (i\partial_t + \Delta_{\mathbb{R}\times\mathbb{T}^3})u = -|u|^2 u, \\ u(0,x) = u_0 \in H^1(\mathbb{R} \times \mathbb{T}^3). \end{cases} \tag{1.1}$$

Generally, the equation (1.1) can be considered as a special case of the general NLS posed on the product spaces $\mathbb{R}^m \times \mathbb{T}^n$, where $m, n \in \mathbb{N}^+$:

$$\begin{cases} (i\partial_t + \Delta_{\mathbb{R}^m\times\mathbb{T}^n})u = \pm|u|^{p-1} u, \\ u(0,x) = u_0 \in H^1(\mathbb{R}^m \times \mathbb{T}^n). \end{cases} \tag{1.2}$$

Here the product spaces $\mathbb{R}^m \times \mathbb{T}^n$ are known as 'semiperiodic spaces' as well as 'waveguide manifolds' (or waveguides for short). The initial value problem (1.2) is called defocusing if the sign of the nonlinearity is positive while it is called focusing if the sign is negative. Throughout this paper, we focus on the focusing, cubic model (1.1). There are three important conserved quantities of (1.1) as follows.

$$\text{Mass:} \quad M(u(t)) = \int_{\mathbb{R}\times\mathbb{T}^3} |u(t,x,y)|^2 \, dxdy, \tag{1.3}$$

$$\text{Energy:} \quad E(u(t)) = \int_{\mathbb{R}\times\mathbb{T}^3} \frac{1}{2}|\nabla u(t,x,y)|^2 - \frac{1}{4}|u(t,x,y)|^4 \, dxdy, \tag{1.4}$$

$$\text{Momentum:} \quad P(u(t)) = \operatorname{Im} \int_{\mathbb{R}\times\mathbb{T}^3} \overline{u(t,x,y)} \nabla u(t,x,y) \, dxdy. \tag{1.5}$$

**Remark 1.1.** The conservation laws are consistent with the analogues in the Euclidean and tori cases. Here we point out that a significant difference between defocusing and focusing models is the sign of the potential energy in (1.4). In the focusing scenario, negative energies are allowed.

The well-posedness theory and the long time behavior of NLS are important matters in the area of dispersive PDEs and have been studied widely in recent decades. Naturally, the Euclidean case is studied first and the theory, at least in the defocusing setting, is well established. We refer to [7, 11, 12, 15, 24, 26, 30, 32] for some important Euclidean results. Moreover, we refer to [19, 21, 27, 33] for works on tori. We may roughly think of the waveguide case as the 'intermediate point between the Euclidean case and the tori case since the waveguide manifold is the product of the Euclidean space and the tori. The techniques used in Euclidean and tori settings are frequently combined and applied to the waveguides problems. We refer to [8, 9, 10, 16, 19, 20, 21, 22, 33, 35, 34, 36] for some NLS results in the waveguide setting.

[1] X.Y. is funded in part by the Jarve Seed Fund and an AMS-Simons travel grant.

Generally, when one considers long time dynamics of NLS with large data, the focusing model is quite different from the defocusing one (although the differences can not been seen in the local theory). Also the dynamics of NLS are much richer in the focusing setting, which complicates the problem. In fact, we could expect many things to happen, such as finite/infinite blow-up and soliton-like behaviors. We refer to [15, 24, 28] regarding some focusing NLS results in Euclidean spaces. However, general long time dynamics for focusing NLS are not well understood in both tori and waveguides. Here we refer to [33] for a focusing result on tori. The analogue in waveguides remains completely open. To the best of authors knowledge, we believe that the current paper is the first result on understanding long time dynamics for the focusing NLS within the context of waveguides. Before stating our main theorem, we also give references for NLS in other manifolds such as hyperbolic spaces and spheres [23, 29] and refer to [6, 14, 26, 32] for some classical textbooks or notes on dispersive PDEs.

We now state the main results of this paper .

**Theorem 1.2** (Global well-posedness result)**.** *Consider the initial value problem* (1.1)*. Assume the solution $u$ to* (1.1) *with the maximal lifespan $I$ satisfies*

$$\sup_{t\in I} \|u(t)\|_{\dot{H}^1(\mathbb{R}\times\mathbb{T}^3)} < \|W\|_{\dot{H}^1(\mathbb{R}^4)}, \tag{1.6}$$

*where $W$ is the ground state on $\mathbb{R}^4$. Then* (1.1) *is globally well-posed on $\mathbb{R}\times\mathbb{T}^3$.*

**Remark 1.3.** Note that (1.6) in Theorem 1.2 is an *a priori* assumption on the solution. In some cases, it is natural and useful to transfer such assumption directly on the initial data. More precisely, we can replace the assumption (1.6) by either (3.17) or (3.18) via an energy trapping argument. See Section 3 for details. The main theorem (Theorem 1.2) together with the energy trapping lemma (Theorem 3.7) implies the following Corollary 1.4.

**Corollary 1.4.** *The initial value problem* (1.1) *is globally well-posed if the initial datum $u_0$ satisfies the assumption* (3.17) *or assumption* (3.18)*.*

**Remark 1.5.** We may compare Theorem 1.2 to Theorem 1.1 in Ionescu-Pausader [22] which concerned the defocusing case. In their case, such an assumption is not required, since the conserved energy in the defocusing setting is positive definite and controls the $\dot{H}^1$ norm of the solution. However, such good control fails in the focusing scenario.

The proof of the Theorem 1.2 follows from a standard skeleton based on the famous Kenig-Merle machinery (see [24, 25] for references), which is well known as 'Concentration compactness/Rigidity method'. In the waveguide setting, our main technical ingredients include adapted atomic function spaces, Strichartz estimates, Euclidean NLS solution approximations and the profile decomposition.

After presenting the main results, let us focus on the similarities and differences between this current result and previous related works. There are three major differences from the defocusing setting (see [22]). In fact, these three differences share a deep reason – the focusing nature.

**(1) An energy trapping argument based on a sharp Sobolev inequality.** As we mentioned above, the energy is no longer positive definite in our case, hence not controlling the $\dot{H}^1$ norm of the solution. It is also known that the negative energy (when data are beyond the the ground state) may lead to blow up solutions. So to best characterize the positive region of the energy, we need to guarantee that the solution stays below the ground state all the time. To this end, we employ a well designed energy trapping argument, that is, up to certain modification on either the energy or the $\dot{H}^1$ norm, the assumption of initial data being below the ground state implies the (modified) energy stays in this good region within its maximal lifespan. Here the modification is necessary and directly connected to the sharp Sobolev inequality. In fact, such Sobolev inequality enters initially when we classify the initial data, thus we have to carry it to the energy trapping argument, hence resulting in the modification. Note that in the defocusing case, all the things we discussed here are not needed at all.

**(2) Euclidean profile approximations.** Due to the lack of control in the $\dot{H}^1$ norm of the solution, our approximations need to carry the *a priori* assumption (1.6) everywhere. Another by-product here is that the linear profile decomposition is not enough to analysis the structure of the nonlinear solutions. This is because the *a priori* assumption is not sole on one particular time, but on the maximal lifespan. Thus we have to develop a corresponding decomposition for the nonlinear profiles (Proposition 5.6) and establish a better understanding on their almost orthogonal property (Lemma 5.9). Notice that this step is not necessary in the defocusing setting.

**(3) The rigidity theorem.** As a consequence of what we discussed above, the main difference in the rigidity argument from the defocusing setting is that we employ the nonlinear profile decomposition and rely on the nonlinear interaction of these profiles, especially the almost orthogonality of the nonlinear profiles.

These differences will be discussed explicitly in the following sections respectively.

We expect that one can obtain analogous focusing results for NLS on other waveguide manifolds using similar methods. One potential difference is the profile decomposition since it is tightly dependent on the geometric structure. Also some other technical issues may appear.

The organization of the rest of this paper is: in Section 2, we discuss the preliminaries; in Section 3, we present a Sobolev-type inequality in the waveguide setting and the corresponding energy trapping-type argument; in Section 4, we review the road map and the main tools for the defocusing case in [22] which also work in our problem; in Section 5, we show the Euclidean profile approximation for the focusing case; in Section 6, we prove the main theorem by the energy induction argument.

## 2. Preliminaries

In this section, we discuss notations, the ground state $W$ and the function spaces that will be used in this paper. These spaces are initially established by S. Herr, D.Tataru, and N. Tzvetkov [19, 20] and have been especially widely applied for NLS problems on tori and waveguides. See [21, 22, 27, 33, 35, 34] for examples. Also see Dodson [13] for a result which applies similar spaces in the Euclidean setting.

### 2.1. Notations.

We write $A \lesssim B$ to say that there is a constant $C$ such that $A \leq CB$. We use $A \simeq B$ when $A \lesssim B \lesssim A$. Particularly, we write $A \lesssim_u B$ to express that $A \leq C(u)B$ for some constant $C(u)$ depending on $u$.

Throughout this paper, we regularly refer to the spacetime norms

$$\|u\|_{L_t^p L_x^q(I_t\times\mathbb{R}\times\mathbb{T}^3)} = \left(\int_{I_t}\left(\int_{\mathbb{R}\times\mathbb{T}^3} |u(x)|^q\,dx\right)^{\frac{p}{q}} dt\right)^{\frac{1}{p}}.$$

Now we turn to Fourier transformation and Littlewood-Paley Theory. We define the Fourier transform on $\mathbb{R}\times\mathbb{T}^3$ as follows:

$$(\mathcal{F}f)(\xi) = \int_{\mathbb{R}\times\mathbb{T}^3} f(x)e^{-ix\cdot\xi}\,dx, \tag{2.1}$$

where $\xi = (\xi_1, \xi_2, \xi_3, \xi_4) \in \mathbb{R} \times \mathbb{Z}^3$. We also note the Fourier inversion formula

$$f(x) = c \sum_{(\xi_2,\xi_3,\xi_4)\in\mathbb{Z}^3} \int_{\xi_1\in\mathbb{R}} (\mathcal{F}f)(\xi) e^{ix\cdot\xi}\, d\xi_1.$$

We define the Schrödinger propagator $e^{it\Delta}$ by

$$(\mathcal{F}e^{it\Delta} f)(\xi) = e^{-it|\xi|^2} (\mathcal{F}f)(\xi).$$

We are now ready to define the Littlewood-Paley projections. First, we fix $\eta_1 : \mathbb{R} \to [0, 1]$, a smooth even function satisfying

$$\eta_1(\xi) = \begin{cases} 1, & |\xi| \leq 1, \\ 0, & |\xi| \geq 2, \end{cases}$$

and $N = 2^j$ a dyadic integer. Let $\eta^4 = \mathbb{R}^4 \to [0, 1]$, $\eta^4(\xi) = \eta_1(\xi_1)\eta_1(\xi_2)\eta_1(\xi_3)\eta_1(\xi_4)$. We define the Littlewood-Paley projectors $P_{\leq N}$ and $P_N$ by

$$\mathcal{F}(P_{\leq N} f)(\xi) := \eta^4(\frac{\xi}{N})\mathcal{F}(f)(\xi), \quad \xi \in (\mathbb{R} \times \mathbb{Z}^3),$$

and

$$P_N f = P_{\leq N} f - P_{\leq \frac{N}{2}} f.$$

For any $a \in (0, \infty)$ we define

$$P_{\leq a} := \sum_{N\leq a} P_N, \quad P_{>a} := \sum_{N>a} P_N.$$

2.2. **Ground state structure on $\mathbb{R}^4$.** In this subsection, we briefly discuss the structure of ground state $W$ which is introduced and studied in [1, 31]. We consider the focusing energy-critical NLS on Euclidean spaces with dimensions $d \geq 3$,

$$\begin{cases} (i\partial_t + \Delta)u = -|u|^{\frac{4}{d-2}} u, & (t, x) \in I \times \mathbb{R}^d \\ u(0, x) = u_0 \in \dot{H}^1(\mathbb{R}^d). \end{cases}$$

Then there is an important radial stationary solution $W$ satisfying the following elliptic equation:

$$\Delta_{\mathbb{R}^d} W = -|W|^{\frac{4}{d-2}} W,$$

and the explicit expression of $W$ is given by

$$W(x) = \frac{1}{(1 + \frac{|x|^2}{d(d-2)})^{\frac{d-2}{2}}}.$$

Ground states are important properties and structures for focusing evolutionary equations. As for our problem, we consider the case $d = 4$ since the whole dimension of the waveguide manifold $\mathbb{R} \times \mathbb{T}^3$ is four. The energy trapping argument related to this ground state will be addressed in Section 3. Moreover, there is an important long time dynamics result in [15] for focusing energy critical NLS which will be discussed and used in Section 5.

2.3. **Function spaces.** In this subsection, we describe the function spaces used in this paper.

**Definition 2.1** ($U^p$ spaces)**.** Let $1 \leq p < \infty$, and $H$ be a complex Hilbert space. A $U^p$-atom is a piecewise defined function, $a : \mathbb{R} \to H$,

$$a = \sum_{k=1}^{K} \chi_{[t_{k-1},t_k)} \phi_{k-1}$$

where $\{t_k\}_{k=0}^K \in \mathcal{Z}$ and $\{\phi_k\}_{k=0}^{K-1} \subset H$ with $\sum_{k=0}^K \|\phi_k\|_H^p = 1$. Here we let $\mathcal{Z}$ be the set of finite partitions $-\infty < t_0 < t_1 < ... < t_K \leq \infty$ of the real line.

The atomic space $U^p(\mathbb{R};H)$ consists of all functions $u:\mathbb{R}\to H$ such that $u=\sum_{j=1}^{\infty}\lambda_j a_j$ for $U^p$-atoms $a_j$, $\{\lambda_j\}\in l^1$, with norm

$$\|u\|_{U^p}:=\inf\{\sum_{j=1}^{\infty}|\lambda_j| : u=\sum_{j=1}^{\infty}\lambda_j a_j, \lambda_j\in\mathbb{C},\quad a_j \text{ are } U^p\text{-atoms}\}.$$

**Definition 2.2** ($V^p$ spaces)**.** Let $1\le p<\infty$, and $H$ be a complex Hilbert space. We define $V^p(\mathbb{R},H)$ as the space of all functions $v:\mathbb{R}\to H$ such that

$$\|u\|_{V^p}:=\sup_{\{t_k\}_{k=0}^{K}\in\mathcal{Z}}(\sum_{k=1}^{K}\|v(t_k)-v(t_{k-1})\|_H^p)^{\frac{1}{p}}\le+\infty,$$

where we use the convention $v(\infty)=0$. Also, we denote the closed subspace of all right-continuous functions $v:\mathbb{R}\to H$ such that $\lim_{t\to-\infty}v(t)=0$ by $V_{rc}^p(\mathbb{R},H)$.

To best serve for this problem, we choose the Hilbert space $H$ to be the $L^2$-based Sobolev spaces $H^s(\mathbb{R}\times\mathbb{T}^3)$.

**Definition 2.3** ($U_\Delta^p$ and $V_\Delta^p$ spaces)**.** For $s\in\mathbb{R}$, we let $U_\Delta^pH^s$ resp. $V_\Delta^pH^s$ be the spaces of all functions such that $e^{-it\Delta}u(t)$ is in $U^p(\mathbb{R},H^s)$ resp. $V_{rc}^p(\mathbb{R},H)$, with norms

$$\|u\|_{U_\Delta^pH^s}=\|e^{-it\Delta}u\|_{U^p(\mathbb{R},H^s)},\quad \|u\|_{V_\Delta^pH^s}=\|e^{-it\Delta}u\|_{V^p(\mathbb{R},H^s)}.$$

For $C=[-\frac{1}{2},\frac{1}{2})^4\in\mathbb{R}^4$ and $z\in\mathbb{R}^4$, we denote by $C_z=z+C$ the translate by $z$ and define the sharp projection operator $P_{C_z}$ as follows: (recall $\mathcal{F}$ is the Fourier transform defined in (2.1)):

$$\mathcal{F}(P_{C_z}f)=\chi_{C_z}(\xi)\mathcal{F}(f)(\xi).$$

Here $\chi_{C_z}$ is the characteristic function restrained on $C_z$. Using the same modifications for the atomic and variation space norms, we can introduce the following intermediate norms.

**Definition 2.4** ($X^s$ and $Y^s$ norms)**.** For $s\in\mathbb{R}$, we define:

$$\|u\|_{X^s(\mathbb{R})}^2=\sum_{z\in\mathbb{Z}^4}\langle z\rangle^{2s}\|P_{C_z}u\|_{U_\Delta^2(\mathbb{R};L^2)}^2$$

and similarly we have,

$$\|u\|_{Y^s(\mathbb{R})}^2=\sum_{z\in\mathbb{Z}^4}\langle z\rangle^{2s}\|P_{C_z}u\|_{V_\Delta^2(\mathbb{R};L^2)}^2.$$

For an interval $I\subset\mathbb{R}$, we can also define the restriction spaces $X^s(I)$ and $Y^s(I)$ in the natural way:

$$X^s(I):=\{u\in C(I:H^s):\|u\|_{X^s}:=\sup_{J\subset I,|J|\le1}\inf\{\|v\|_{X^s(\mathbb{R})}:v|_J=u\}\}.$$

Then $Y^s(I)$ norm can be defined in a similar manner.

In fact, both $X^s$ and $Y^s$ norms are stronger than the $L^\infty(\mathbb{R};H^s)$ norm and weaker than the norm $U_\Delta^2(\mathbb{R}:H^s)$. Moreover, we have the following embedding.

**Proposition 2.5** (Embedding [19, 20])**.** *For $p>2$,*

$$U_\Delta^2(\mathbb{R}:H^s)\hookrightarrow X^s\hookrightarrow Y^s\hookrightarrow V_\Delta^2(\mathbb{R}:H^s)\hookrightarrow U_\Delta^p(\mathbb{R}:H^s)\hookrightarrow L^\infty(\mathbb{R};H^s).$$

In order to deal with the nonlinearity, we bring in the following $N$-norm for nonlinear estimates.

**Definition 2.6** ($N$-norm)**.** On a time interval $I$, we define

$$\|h\|_{N^s(I)}=\|\int_a^t e^{i(t-s)\Delta}h(s)\,ds\|_{X^s(I)}.$$

The last norm that we need in this paper is the solution norm, i.e. '$Z$-norm'.

**Definition 2.7** ($Z$-norm)**.** On a time interval $I$, we define

$$\|u\|_{Z(I)} = \sup_{J\subset I, |J|\leq 1} (\sum_{N\geq 1} N^2 \|P_N u\|^4_{L^4_{t,x,y}(J\times\mathbb{R}\times\mathbb{T}^3)})^{\frac{1}{4}}.$$

We note that $Z$ is a weaker norm than $X^1$ in in following sense,

$$\|u\|_{Z(I)} \lesssim \|u\|_{X^1(I)},$$

which follows from Strichartz estimates (see Section 4).

We close this section by presenting the following duality property.

**Proposition 2.8** (Proposition 2.11 in [19])**.** *If $f \in L^1_t(I, H^1(\mathbb{R}\times\mathbb{T}^3))$, then*

$$\|f\|_{N(I)} \lesssim \sup_{v\in Y^{-1}(I), \|v\|_{Y^{-1}(I)}\leq 1} \int_{I\times(\mathbb{R}\times\mathbb{T}^3)} f(t,x)\overline{v(t,x)}\, dxdt.$$

*Also, we have the following estimate holds for any smooth function $g$ on an interval $I = [a,b]$:*

$$\|g\|_{X^1(I)} \lesssim \|g(a)\|_{H^1(\mathbb{R}\times\mathbb{T}^3)} + (\sum_N \|P_N(i\partial_t + \Delta)g\|^2_{L^1_t(I,H^1(\mathbb{R}\times\mathbb{T}^3))})^{\frac{1}{2}}.$$

Note that for simplicity, all the spacetime norms that we will use in the rest of this paper will be restricted within the time interval $|I| \leq 1$.

## 3. Energy trapping argument

In this section, we aim to impose suitable conditions on the initial data in (1.1) to obtain the assumption (1.6) in Theorem 1.2. Similar as in Theorem 3.9 of [24], we will prove a energy trapping argument which is based on a sharp Sobolev inequality in the setting of waveguide manifold.

### 3.1. Sharp Sobolev inequality for waveguide manifolds.

**Proposition 3.1.** *For a function in the waveguide $\mathbb{R}^n \times \mathbb{T}^{d-n}$ ($1 \leq n < d$) with whole dimension $d \geq 2$, there exists a constant $c > 0$ such that*

$$\|f(x,y)\|_{L^{p^*}(\mathbb{R}^n\times\mathbb{T}^{d-n})} \leq C_d\|\nabla f\|_{L^2} + c\|f\|_{L^2}, \tag{3.1}$$

*where $p^* = \frac{2d}{d-2}$ and $C_d$ is the best constant for sharp Sobolev inequality in the Euclidean setting. Moreover the constant $C_d$ is sharp.*

**Remark 3.2.** The Sobolev inequality for waveguide manifolds is almost identical to the tori analogue when replacing the $d$-dimensional tori by a waveguide with the same whole dimension $d$. See [17, 18] for the tori case.

**Remark 3.3.** It is worth pointing out that Proposition 3.1 is sharp in term of the constant $C_d$, while $c$ in the second term is not optimal. For the existence of such $c$, obviously we can take the Sobolev constant $c > 0$ in its corresponding Euclidean case.

*Proof of Proposition 3.1.* First, Proposition 3.1 is a direct consequence of the sharp Sobolev inequality in the Euclidean setting by taking a spatial truncation function. Thus it suffices to show the sharpness of (3.1).

Assume that there exists a smaller constant $0 < K_d < C_d$ satisfying (3.1) in the following sense,

$$\|f(x,y)\|_{L^{p^*}(\mathbb{R}^n\times\mathbb{T}^{d-n})} \leq K_d\|\nabla f\|_{L^2} + c\|f\|_{L^2}.$$

Then we take a smooth function $f$ which is compactly supported in a small enough ball $B_\delta$ (with radius $\delta$ to be decided later). Moreover, on $B_\delta$, using Hölder inequality, we have

$$\|f\|_{L^2(B_\delta)} \le |B_\delta|^{\frac{d+2}{4d}} \cdot \|f\|_{L^{p^\star}(B_\delta)}.$$

By choosing a small enough $\delta$, the $L^2$ term on the right hand side can be absorbed by the term on the left. Thus we can choose $K_d^{'}$ satisfying $K_d < K_d^{'} < C_d$ such that for function $f$ supported on $B_\delta$,

$$\|f\|_{L^{p^*}(\mathbb{R}^d)} \le K_d^{'} \|\nabla f\|_{L^2(\mathbb{R}^d)}.$$

Noe we rescale it back to the general compactly supported function in $\mathbb{R}^d$. Let $f \in C_0^\infty(\mathbb{R}^d)$ and set $u_\lambda(x) = f(\lambda x)$ $(\lambda > 0)$, then for $\lambda$ large enough, $f_\lambda \in C_0^\infty(B_\delta)$. Hence,

$$\|f_\lambda\|_{L^{p^*}(\mathbb{R}^d)} \le K_d^{'} \|\nabla f_\lambda\|_{L^2(\mathbb{R}^d)}. \tag{3.2}$$

By rescaling back, (3.2) implies that

$$\|f\|_{L^{p^*}(\mathbb{R}^d)} \le K_d^{'} \|\nabla f\|_{L^2(\mathbb{R}^d)},$$

which contradicts the sharpness of the Sobolev inequality in the Euclidean setting.

The proof of Proposition 3.1 is complete. □

Now we consider the estimate for the specific case, i.e. $\mathbb{R}\times\mathbb{T}^3$ as follows.

**Corollary 3.4** (Sobolev embedding in $\mathbb{R}\times\mathbb{T}^3$). *Let $f \in H^1(\mathbb{R}\times\mathbb{T}^3)$, then there exists a positive constant $c_*$, such that*

$$\|f\|^2_{L^4(\mathbb{R}\times\mathbb{T}^3)} \le C_4^2(\|f\|^2_{\dot H^1(\mathbb{R}\times\mathbb{T}^3)} + c_*\|f\|^2_{L^2(\mathbb{R}\times\mathbb{T}^3)}).$$

*where $C_4$ is the optimal in this inequality.*

**Remark 3.5.** Note that $\|u\|^2_{H^1_*(\mathbb{R}\times\mathbb{T}^3)} = \|u\|^2_{\dot H^1(\mathbb{R}\times\mathbb{T}^3)} + c_*\|u\|^2_{L^2(\mathbb{R}\times\mathbb{T}^3)}$, then the Sobolev embedding (Corollary 3.4) can be also written in the form:

$$\|f\|_{L^4(\mathbb{R}\times\mathbb{T}^3)} \le C_4\|f\|_{H^1_*(\mathbb{R}\times\mathbb{T}^3)}.$$

Also $C_4$ can calculating from $W$ as follows:

$$\|W\|^2_{\dot H^1(\mathbb{R}^4)} = \|W\|^4_{L^4(\mathbb{R}^4)} := \frac{1}{C_4^4} \quad \text{and then} \quad E_{\mathbb{R}^4}(W) = \frac{1}{4C_4^4}. \tag{3.3}$$

3.2. **An energy trapping lemma.** Given this different Sobolev embedding, in order to run a suitable energy trapping argument in the focusing case, let us introduce two modified energies of $u$:

$$E_*(u)(t) := \frac{1}{2}(\|u(t)\|^2_{\dot H^1(\mathbb{R}\times\mathbb{T}^3)} + c_*\|u(t)\|^2_{L^2(\mathbb{R}\times\mathbb{T}^3)}) - \frac{1}{4}\|u(t)\|^4_{L^4(\mathbb{R}\times\mathbb{T}^3)}, \tag{3.4}$$

and

$$E_{**}(u)(t) := \frac{1}{2}(\|u(t)\|^2_{\dot H^1(\mathbb{R}\times\mathbb{T}^3)} + c_*\|u(t)\|^2_{L^2(\mathbb{R}\times\mathbb{T}^3)}) - \frac{1}{4}\|u(t)\|^4_{L^4(\mathbb{R}\times\mathbb{T}^3)} + \frac{c_*^2C_4^4}{4}\|u(t)\|^4_{L^2(\mathbb{R}\times\mathbb{T}^3)}, \tag{3.5}$$

where $c_*$ is a fixed constant in Corollary 3.4. Notice that by the definitions (3.4), (3.5) and conservation laws, both $E_*(u)(t)$ and $E_{**}(u)(t)$ are conserved in time.

Moreover, we have the following properties of the modified energy defined above.

**Lemma 3.6.** *Suppose $f \in H^1(\mathbb{R}\times\mathbb{T}^3)$.*

**Part (i)** *Assume*

$$\|f\|_{H^1_*(\mathbb{R}\times\mathbb{T}^3)} < \|W\|_{\dot H^1(\mathbb{R}^4)} \quad \text{and} \quad E_*(f) < (1-\delta_0)E_{\mathbb{R}^4}(W), \tag{3.6}$$

*for $\delta_0 > 0$, then there exists $\bar\delta = \bar\delta(\delta_0) > 0$ such that*

$$\|f\|^2_{H^1_*(\mathbb{R}\times\mathbb{T}^3)} < (1-\bar\delta)\|W\|^2_{\dot H^1(\mathbb{R}^4)} \tag{3.7}$$

$$\|f\|^2_{H^1_*(\mathbb{R}\times\mathbb{T}^3)} - \|f\|^4_{L^4(\mathbb{R}\times\mathbb{T}^3)} \geq \bar{\delta}\|f\|^2_{H^1_*(\mathbb{R}\times\mathbb{T}^3)}, \tag{3.8}$$

*and in particular*

$$E_*(f) \geq \frac{1}{4}(1+\bar{\delta})\|f\|^2_{H^1_*(\mathbb{R}\times\mathbb{T}^3)}. \tag{3.9}$$

**Part (ii)** *Assume*

$$\|f\|_{\dot{H}^1(\mathbb{R}\times\mathbb{T}^3)} < \|W\|_{\dot{H}^1(\mathbb{R}^4)} \qquad \text{and} \qquad E_{**}(f) < (1-\delta_0)E_{\mathbb{R}^4}(W), \tag{3.10}$$

*for $\delta_0 > 0$, then there exists $\bar{\delta} = \bar{\delta}(\delta_0) > 0$ such that*

$$\|f\|^2_{\dot{H}^1(\mathbb{R}\times\mathbb{T}^3)} < (1-\bar{\delta})\|W\|^2_{\dot{H}^1(\mathbb{R}^4)} \tag{3.11}$$

$$\|f\|^2_{\dot{H}^1(\mathbb{R}\times\mathbb{T}^3)} - \|f\|^4_{L^4(\mathbb{R}\times\mathbb{T}^3)} + 2c_*\|f\|_{L^2(\mathbb{R}\times\mathbb{T}^3)} + c_*^2C_4^4\|f\|^4_{L^2(\mathbb{R}\times\mathbb{T}^3)} \geq \bar{\delta}\|f\|^2_{\dot{H}^1(\mathbb{R}\times\mathbb{T}^3)}, \tag{3.12}$$

*and in particular*

$$E_{**}(f) \geq \frac{1}{4}(1+\bar{\delta})\|f\|^2_{\dot{H}^1(\mathbb{R}\times\mathbb{T}^3)}. \tag{3.13}$$

*Proof of Lemma 3.6.* The proofs of both **Part (i)** and **Part (ii)** are similar, so we will only present the proof of **Part (i)**. First, we follow the proof of Lemma 3.4 in [24], but use the modified Sobolev norm $H^1_*(\mathbb{R}\times\mathbb{T}^3)$ norm instead of the regular $\dot{H}^1(\mathbb{R}\times\mathbb{T}^3)$ norm. Let us first set up an adapted quadratic function $g_1 = \frac{1}{2}y - \frac{C_4^4}{4}y^2$ and plug $y = \|f\|^2_{H^1_*(\mathbb{R}\times\mathbb{T}^3)}$ in it. By Corollary 3.4 and the assumption (3.6), we have then that

$$\begin{aligned} g_1(\|f\|^2_{H^1_*}) &= \frac{1}{2}\|f\|^2_{H^1_*} - \frac{C_4^4}{4}\|f\|^4_{H^1_*} \\ &\leq \frac{1}{2}\|f\|^2_{H^1_*} - \frac{1}{4}\|f\|^4_{L^4} = E_*(f) \\ &< (1-\delta_0)E_{\mathbb{R}^4}(W) = (1-\delta_0)g_1(\|W\|^2_{\dot{H}^1(\mathbb{R}^4)}) \end{aligned} \tag{3.14}$$

which applies (3.7) together with the quadratic function $g_1$ when $\bar{\delta} \sim \delta_0^{\frac{1}{2}}$.

Then we set up another quadratic function $g_2(y) = y - C_4^4y^2$. After plugging in $y = \|f\|^2_{H^1_*(\mathbb{R}\times\mathbb{T}^3)}$ into $g_2$, by Corollary 3.4 again, we have that

$$g_2(\|f\|^2_{H^1_*(\mathbb{R}\times\mathbb{T}^3)}) = \|f\|^2_{H^1_*(\mathbb{R}\times\mathbb{T}^3)} - C_4^4\|f\|^4_{H^1_*(\mathbb{R}\times\mathbb{T}^3)} \leq \|f\|^2_{H^1_*(\mathbb{R}\times\mathbb{T}^3)} - \|f\|^4_{L^4(\mathbb{R}\times\mathbb{T}^3)}. \tag{3.15}$$

It is easy to check that $g_2(0) = 0$, $g_2''(y) = -2C_4^4 < 0$ and $\|f\|^2_{H^1_*(\mathbb{R}\times\mathbb{T}^3)} < (1-\bar{\delta})\|W\|^2_{\dot{H}^1(\mathbb{R}^4)}$. By Jensen's inequality and (3.3), we have

$$g_2(\|f\|^2_{H^1_*(\mathbb{R}\times\mathbb{T}^3)}) > g_2((1-\bar{\delta})\|W\|^2_{\dot{H}^1(\mathbb{R}^4)})\frac{\|f\|^2_{H^1_*(\mathbb{R}\times\mathbb{T}^3)}}{(1-\bar{\delta})\|W\|^2_{\dot{H}^1(\mathbb{R}^4)}} = \bar{\delta}\|f\|^2_{H^1_*(\mathbb{R}\times\mathbb{T}^3)}. \tag{3.16}$$

Based on (3.15) and (3.16), we prove (3.8), and then (3.9) follows

$$E_*(f) = \frac{1}{4}\|f\|^2_{H^1_*(\mathbb{R}\times\mathbb{T}^3)} + \frac{1}{4}(\|f\|^2_{H^1_*(\mathbb{R}\times\mathbb{T}^3)} - \|f\|^4_{L^4(\mathbb{R}\times\mathbb{T}^3)}) \geq \frac{1}{4}(1+\bar{\delta})\|f\|^2_{H^1_*(\mathbb{R}\times\mathbb{T}^3)}.$$

□

We now present the energy trapping arguments corresponding to the two parts in Lemma 3.6.

**Theorem 3.7** (Energy trapping)**.** *Let $u$ be a solution of IVP* (1.1).

**Part (i)** *Suppose*

$$\|u_0\|_{H^1_*(\mathbb{R}\times\mathbb{T}^3)} < \|W\|_{\dot{H}^1(\mathbb{R}^4)}, \quad E_*(u_0) < (1-\delta_0)E_{\mathbb{R}^4}(W); \tag{3.17}$$

*for some* $\delta_0 > 0$. *Let* $I \ni 0$ *be the maximal interval of existence, then there exists* $\bar{\delta} = \bar{\delta}(\delta_0) > 0$ *such that for all* $t \in I$

$$\|u(t)\|^2_{H^1_*(\mathbb{R}\times\mathbb{T}^3)} < (1-\bar{\delta})\|W\|_{\dot{H}^1(\mathbb{R}^4)},$$
$$\|u(t)\|^2_{H^1_*(\mathbb{R}\times\mathbb{T}^3)} - \|u(t)\|^4_{L^4(\mathbb{R}\times\mathbb{T}^3)} \geq \bar{\delta}\|u(t)\|^2_{H^1_*(\mathbb{R}\times\mathbb{T}^3)},$$

*and in particular*

$$E_*(u)(t) \geq \frac{1}{4}(1+\bar{\delta})\|u(t)\|^2_{H^1_*(\mathbb{R}\times\mathbb{T}^3)}.$$

**Part (ii)** *Suppose*

$$\|u_0\|_{\dot{H}^1(\mathbb{R}\times\mathbb{T}^3)} < \|W\|_{\dot{H}^1(\mathbb{R}^4)}, \quad E_{**}(u_0) < (1-\delta_0)E_{\mathbb{R}^4}(W); \tag{3.18}$$

*for some* $\delta_0 > 0$. *Let* $I \ni 0$ *be the maximal interval of existence, then there exists* $\bar{\delta} = \bar{\delta}(\delta_0) > 0$ *such that for all* $t \in I$

$$\|u(t)\|^2_{\dot{H}^1(\mathbb{R}\times\mathbb{T}^3)} < (1-\bar{\delta})\|W\|_{\dot{H}^1(\mathbb{R}^4)},$$
$$\|u(t)\|^2_{\dot{H}^1(\mathbb{R}\times\mathbb{T}^3)} - \|u(t)\|^4_{L^4(\mathbb{R}\times\mathbb{T}^3)} + 2c_*\|u(t)\|^2_{L^2(\mathbb{R}\times\mathbb{T}^3)} + c_*^2C_4^4\|u(t)\|^4_{L^2(\mathbb{R}\times\mathbb{T}^3)} \geq \bar{\delta}\|u(t)\|^2_{\dot{H}^1(\mathbb{R}\times\mathbb{T}^3)},$$

*and in particular*

$$E_{**}(u)(t) \geq \frac{1}{4}(1+\bar{\delta})\|u(t)\|^2_{\dot{H}^1(\mathbb{R}\times\mathbb{T}^3)}.$$

*Proof of Theorem 3.7.* By the conservation laws of modified energies $E_*(u)(t)$ and $E_{**}(u)(t)$, this theorem is proven directly by Lemma 3.6 and a continuity argument. □

In summary, Theorem 3.7 shows that if the initial datum satisfies the condition (3.17) or (3.18) then the solution $u(t)$ satisfies $\|u(t)\|_{\dot{H}^1(\mathbb{R}\times\mathbb{T}^3)} < \|W\|_{\dot{H}^1(\mathbb{R}^4)}$ for all $t$ in the lifespan of the solution. So Corollary 1.4 holds due to the same reasoning in Theorem 1.2.

## 4. An overview of results for the defocusing analogue

In this section, we review some important theorems and the road map of the proof for the defocusing case as in [22]. They also play a significant role in the focusing setting.

### 4.1. Strichartz Estimate.

**Proposition 4.1** (Strichartz estimate). *For any* $p > p_0 = 3$, $N \geq 1$ *and* $f \in L^2(\mathbb{R}\times\mathbb{T}^3)$,

$$\|e^{it\Delta}P_N f\|_{L^p_{t,x}([-1,1]\times\mathbb{R}\times\mathbb{T}^3)} \lesssim N^{2-\frac{6}{p}}\|f\|_{L^2(\mathbb{R}\times\mathbb{T}^3)}.$$

**Remark 4.2.** In [22], the threshold $p_0$ is $\frac{18}{5}$. Later, it is improved to 3 in [2] based on Bourgain-Demeter's [5] decoupling method. See also [3, 4, 27] for the tori analogue.

4.2. **Local Theory.** Based on the Strichartz estimate and the properties of function spaces, useful nonlinear estimates can be established. Furthermore, we can obtain the local well-posedness theory, which works for both of defocusing and focusing cases. We recall the local theory as follows. (Nonlinear estimates are dismissed. The proofs are contained in Section 3 of [22].)

**Theorem 4.3** (Local well-posedness). *Let $E > 0$ and $\|u_0\|_{H^1(\mathbb{R}\times\mathbb{T}^3)} < E$, then there exists $\delta_0 = \delta_0(E) > 0$ such that if*

$$\|e^{it\Delta}u_0\|_{Z(I)}^{\frac{3}{4}} \cdot \|e^{it\Delta}u_0\|_{X^1(I)}^{\frac{1}{4}} < \delta$$

*for some $\delta \leq \delta_0$ and some interval $I \ni 0$ with $|I| \leq 1$, then there exists a unique strong solution $u \in X^1(I)$ with initial datum $u(0) = u_0$ satisfying*

$$\|u(t) - e^{it\Delta_{\mathbb{R}\times\mathbb{T}^3}}u_0\|_{X^1(I)} \leq \delta^{\frac{5}{3}}.$$

**Theorem 4.4** (Controlling norm). *Let $u \in X^1(I)$ be a strong solution on a bounded open interval $I \in \mathbb{R}$ satisfying*

$$\|u\|_{Z(I)} < \infty.$$

*Then we have that there exists an open interval $J$ with $I \subset J$ such that $u$ can be extended to a strong solution of* (1.1) *on $J$. In particular, if $u$ blows up in finite time, then $u$ blows up in the $Z$-norm.*

**Theorem 4.5** (Stability theory). *Let $I \in \mathbb{R}$ be a bounded open interval, and let $\widetilde{u} \in X^1(I)$ solve the approximate equation,*

$$(i\partial_t + \Delta_{\mathbb{R}\times\mathbb{T}^3})\widetilde{u} = \rho\,|\widetilde{u}|^2\widetilde{u} + e \text{ where } \rho \in \{-1, 0, 1\}.$$

*Assume in addition that,*

$$\|\widetilde{u}\|_{Z(I)} + \|\widetilde{u}\|_{L_t^\infty(I,H^1(\mathbb{R}\times\mathbb{T}^3))} \leq M,$$

*for some $M \in [1, +\infty)$. Assume $t_0 \in I$ and $u_0 \in H^1(\mathbb{R}\times\mathbb{T}^3)$ is such that the smallness condition*

$$\|\widetilde{u}(t_0) - u_0\|_{H^1(\mathbb{R}\times\mathbb{T}^3)} + \|e\|_{N(I)} \leq \varepsilon$$

*holds for some $0 < \varepsilon < \varepsilon_1$, where $\varepsilon_1 \leq 1$ is a small constant $\varepsilon_1 = \varepsilon_1(M) > 0$. Then there exists a solution $u(t)$ to the exact equation:*

$$(i\partial_t + \Delta_{\mathbb{R}\times\mathbb{T}^3})u = \rho\,|u|^2u$$

*with initial datum $u(t_0) = u_0$ such that*

$$\|u\|_{X^1(I)} + \|\widetilde{u}\|_{X^1(I)} \leq C(M), \qquad \|u - \widetilde{u}\|_{X^1(I)} \leq C(M)\varepsilon.$$

4.3. **Profile Decomposition.** One of the powerful tools to extend the local theory to global is the 'profile decomposition'. The types of profiles greatly depend on the problem. For our model, there are two types of profiles, i.e. Euclidean profiles and scale-one profiles.

Before introducing Euclidean profiles and scale-one profiles, let us first introduce some related notations. Suppose $f(t,x)$ with $(t,x) \in \mathbb{R}\times\mathbb{T}^3$, $t_0 \in \mathbb{R}$ and $x_0 \in \mathbb{R}\times\mathbb{T}^3$. We define the following two transformations of $f$

$$\begin{aligned} (\pi_{x_0}f)(x) &:= f(x - x_0), \\ (\Pi_{t_0,x_0})f(x) &:= (\pi_{x_0}e^{-it_0\Delta}f)(x). \end{aligned} \tag{4.1}$$

We fix a spherically symmetric function $\eta \in C_0^\infty(\mathbb{R}^4)$ supported in the ball of radius 2 and equal to 1 in the ball of radius 1. Given $\phi \in \dot{H}^1(\mathbb{R}^4)$ and $N \geq 1$, we define

$$T_N\phi(x) := N(Q_N\phi)(N\Psi^{-1}(x)), \text{ where } (Q_N\phi)(y) := \eta(y/N^{\frac{1}{2}})\phi(y), \tag{4.2}$$

where $\Psi : \{x \in \mathbb{R}^4 : |x| < 1\} \to O_0 \subseteq \mathbb{R}\times\mathbb{T}^3, \Psi(x) = x$. And we then have that $T_N : \dot{H}^1(\mathbb{R}^4) \to H^1(\mathbb{R}^4)$ is a linear operator with $\|T_N\phi\|_{H^1(\mathbb{R}\times\mathbb{T}^3)} \lesssim \|\phi\|_{\dot{H}^1(\mathbb{R}^4)}$.

**Definition 4.6** (Frames and Profiles). **Part (i)** We define a frame to be sequence $(N_k, t_k, x_k)_k \in 2^{\mathbb{Z}} \times \mathbb{R} \times (\mathbb{R}\times\mathbb{T}^3)$ . And we can define some types of profiles as follows.

- **(a)** A Euclidean frame is a sequence $\mathcal{F}_e = (N_k, t_k, x_k)$ with $N_k \geq 1, N_k \to \infty, t_k \in \mathbb{R}, x_k \in \mathbb{R}\times\mathbb{T}^3$.
- **(b)** A scale-one frame is a sequence $\mathcal{F}_1 = (1, t_k, x_k)$ with $t_k \in \mathbb{R}, x_k \in \mathbb{R}\times\mathbb{T}^3$.

**Part (ii)** We say that two frames $(N_k,t_k,x_k)_k$ and $(M_k,s_k,y_k)_k$ are orthogonal if

$$\lim_{k\to+\infty}(|\ln\frac{N_k}{M_k}|+N_k^2|t_k-s_k|+N_k|x_k-y_k|)=+\infty.$$

**Part (iii)** We associate a profile defined as:

**(a)** If $\mathcal{O}=(N_k,t_k,x_k)_k$ is a Euclidean frame and for $\varphi\in\dot{H}^1(\mathbb{R}^4)$ we define the Euclidean profile associated to $(\varphi,\mathcal{O})$ as the sequence $\widetilde{\varphi}_{\mathcal{O},k}$ with

$$\widetilde{\varphi}_{\mathcal{O},k}=\Pi_{t_k,x_k}(T_{N_k})(x,y).$$

**(b)** If $\mathcal{O}=(1,t_k,x_k)_k$ is a scale one frame, if $\omega\in H^1(\mathbb{R}\times\mathbb{T}^3)$, we define the scale one profile associated to $(\omega,\mathcal{O})$ as $\widetilde{\omega}_{\mathcal{O},k}$ with

$$\widetilde{\omega}_{\mathcal{O},k}=\Pi_{t_k,x_k}\omega.$$

**Part (iv)** Finally, we say that a sequence of functions $\{f_k\}_k\subset H^1(\mathbb{R}\times\mathbb{T}^3)$ is absent from a frame $\mathcal{O}$ if, up to a subsequence, for any profile $\widetilde{\psi}_{\mathcal{O},k}$ associated with $\mathcal{O}$.

$$\langle f_k,\widetilde{\psi}_{\mathcal{O},k}\rangle_{H^1\times H^1}\to 0\quad\text{as } k\to\infty.$$

**Theorem 4.7** (Profile decomposition). *Assume $\{\phi_k\}_k$ is a sequence of functions satisfying*

$$\|\phi_k\|_{H^1(\mathbb{R}\times\mathbb{T}^3)}<E,$$

*then up to a subsequence, there exists a sequence of Euclidean profiles $\widetilde{\varphi}^{\alpha}_{\mathcal{O}^\alpha,k}$, and scale-one profiles $\widetilde{\omega}^{\beta}_{\mathcal{O}^\beta,k}$ such that, for any $J\geq 0$*

$$\phi_k(x,y)=\sum_{1\leq\alpha\leq J}\widetilde{\varphi}^{\alpha}_{\mathcal{O}^\alpha,k}+\sum_{1\leq\beta\leq J}\widetilde{\omega}^{\beta}_{\mathcal{O}^\beta,k}+R_k^J$$

*where $R_k^J$ is absent from the frames $\mathcal{O}^\alpha$ and satisfies*

$$\limsup_{J\to\infty}\Lambda_\infty(\{R_k^J\})=0.$$

*Here $\Lambda_\infty$ is a functional defined as follows*

$$\Lambda_\infty(\{f_k\})=\limsup_{k\to\infty}\sup_N N^{-1}\|e^{it\Delta}P_Nf_k\|_{L^\infty_{t,x}}.$$

*Additionally, we have the following orthogonality relations*

$$\|\phi_k\|_{L^2}^2=\sum_\alpha\|\widetilde{\varphi}^{\alpha}_{\mathcal{O}^\alpha,k}\|_{L^2}^2+\sum_\beta\|\widetilde{\omega}^{\beta}_{\mathcal{O}^\beta,k}\|_{L^2}^2+\|R_k^J\|_{L^2}^2+o_k(1),$$

$$\|\nabla\phi_k\|_{L^2}^2=\sum_\alpha\|\nabla\widetilde{\varphi}^{\alpha}_{\mathcal{O}^\alpha,k}\|_{L^2}^2+\sum_\beta\|\nabla\widetilde{\omega}^{\beta}_{\mathcal{O}^\beta,k}\|_{L^2}^2+\|\nabla R_k^J\|_{L^2}^2+o_k(1),$$

$$\|\phi_k\|_{L^4}^4=\sum_\alpha\|\widetilde{\varphi}^{\alpha}_{\mathcal{O}^\alpha,k}\|_{L^4}^4+\sum_\beta\|\widetilde{\omega}^{\beta}_{\mathcal{O}^\beta,k}\|_{L^4}^4+o_{J,k}(1).$$

*Recall that our spacetime norms are restricted within the time interval $|I|\leq 1$.*

**Remark 4.8.** We note that, if one regards $\Lambda_\infty$ as a norm, then $Z$-norm is 'weaker' than $\Lambda_\infty$ in the sense that

$$\limsup_{k\to\infty}||f_k||_Z\lesssim\Lambda_\infty(\{f_k\}).$$

We refer to (6.6) in Proposition 6.1 in [22] for the proof. Thus in the above profile decomposition, we have

$$\limsup_{J\to\infty}\limsup_{k\to\infty}\|e^{it\Delta}R_k^J\|_Z=0.$$

Theorem 4.7 is in fact a linear analysis result, thus it works for both defocusing and focusing cases. So far we have reviewed the 'common theories' that both defocusing and focusing settings share, which are originally established in [22]. It is clear that the local theory for the focusing case can be established in a similar way. In the next two sections, we will work on the global argument and then prove the main theorem of this paper.

## 5. Euclidean Profile Approximation

As we explained in the introduction, one of the main differences from the defocusing case is the Euclidean profile approximation. In this section, we will explore the Euclidean profile approximations in the focusing setting.

Let us start by recalling the Euclidean result by Dodson [15].

**Theorem 5.1.** *Assume $\phi \in \dot{H}^1(\mathbb{R}^4)$, and the solution $v(t)$ satisfies*

$$\sup_{t\in I} \|v(t)\|_{\dot{H}^1} < \|W\|_{\dot{H}^1},$$

*where $I$ is the maximal lifespan of the solution. Then there is a unique global solution $v \in C(\mathbb{R} : \dot{H}^1(\mathbb{R}^4))$ of the initial-value problem*

$$(i\partial_t + \Delta_{\mathbb{R}^4})v = -v|v|^2, \quad v(0) = \phi, \tag{5.1}$$

*and*

$$\|\nabla_{\mathbb{R}^4} v\|_{(L_t^\infty L_x^2 \cap L_t^2 L_x^4)(\mathbb{R}\times\mathbb{R}^4)} \leq C(E_{\mathbb{R}^4}(\phi)).$$

*Moreover this solution scatters in the sense that there exists $\phi^{\pm\infty} \in \dot{H}^1(\mathbb{R}^4)$ such that*

$$\|v(t) - e^{it\Delta}\phi^{\pm\infty}\|_{\dot{H}^1(\mathbb{R}^4)} \to 0 \tag{5.2}$$

*as $t \to \pm\infty$. Besides if $\phi \in H^5(\mathbb{R}^4)$, then $v \in C(\mathbb{R} : H^5(\mathbb{R}^4))$ and*

$$\sup_{t\in\mathbb{R}} \|v(t)\|_{H^5(\mathbb{R}^4)} \lesssim_{\|\phi\|_{H^5(\mathbb{R}^4)}} 1.$$

Theorem 5.1 will be used to approximate the nonlinear Euclidean profiles based on the perturbation theory (Theorem 4.5). The propositions and lemmas presented in the rest of this section and their proofs are adapted from those in Section 5 of [22].

Recall that $\eta \in C_0^\infty(\mathbb{R}^4)$ is a spherically symmetric function supported in the ball of radius 2 and equal to 1 in the ball of radius 1. Given $\phi \in \dot{H}^1(\mathbb{R}^4)$ and a real number $N \geq 1$ we define

$$\begin{aligned} Q_N\phi &\in H^1(\mathbb{R}^4), & (Q_N\phi)(x) &= \eta(x/N^{\frac{1}{2}})\phi(x), \\ \phi_N &\in H^1(\mathbb{R}^4), & \phi_N(x) &= N(Q_N\phi)(Nx), \\ f_N &\in H^1(\mathbb{R}\times\mathbb{T}^3), & f_N(y) &= \phi_N(\Psi^{-1}(y)), \end{aligned} \tag{5.3}$$

where $\Psi : \{x \in \mathbb{R}^4 : |x| < 1\} \to O_0 \subseteq \mathbb{R}\times\mathbb{T}^3$, $\Psi(x) = x$. Note that here $f_N$ is just $T_N\phi$.

**Lemma 5.2.** *Assume $\phi \in \dot{H}^1(\mathbb{R}^4)$, $T_0 \in (0,\infty)$, and $\rho \in \{-1, 0\}$ are given, and we define $f_N$ as in (5.3). Under the assumption that for any solutions $v$ of (5.1) with initial datum $\phi$, we have*

$$\sup_{t\in \text{lifespan of } v} \|v(t)\|_{\dot{H}^1(\mathbb{R}^4)} < \|W\|_{\dot{H}^1(\mathbb{R}^4)}.$$

*Then the following conclusions hold:*

**Part (i)** *There is $N_0 = N_0(\phi, T_0)$ sufficiently large such that for any $N \geq N_0$, there is a unique solution $U_N \in C((-T_0N^{-2}, T_0N^{-2}); H^1(\mathbb{R}\times\mathbb{T}^3))$ of the initial-value problem*

$$(i\partial_t + \Delta_{\mathbb{R}\times\mathbb{T}^3})U_N = \rho U_N|U_N|^2, \quad \text{and} \quad U_N(0) = f_N.$$

*Moreover, for any $N \geq N_0$,*

$$\|U_N\|_{X^1(-T_0N^{-2}, T_0N^{-2})} \lesssim_{E_{\mathbb{R}^4}(\phi)} 1.$$

**Part (ii)** *Assume $\varepsilon_1 \in (0, 1]$ is sufficiently small (depending only on $E_{\mathbb{R}^4}(\phi)$), $\phi' \in H^5(\mathbb{R}^4)$, and $\|\phi - \phi'\|_{\dot{H}^1(\mathbb{R}^4)} \leq \varepsilon_1$. Let $v' \in C(\mathbb{R} : H^5)$ denote the solution of the initial-value problem*

$$(i\partial_t + \Delta_{\mathbb{R}^4})v' = \rho v'|v'|^2, \quad v'(0) = \phi',$$

*where* (5.1) *is satisfied when* $\rho=-1$. *For* $R\geq 1$ *and* $N\geq 10R$, *we define*

$$\begin{aligned} v'_R(t,x) &= \eta(x/R)v'(t,x) && (t,x)\in(-T_0,T_0)\times\mathbb{R}^4,\\ v'_{R,N}(t,x) &= Nv'_R(N^2t,Nx) && (t,x)\in(-T_0N^{-2},T_0N^{-2})\times\mathbb{R}^4,\\ V_{R,N}(t,y) &= v'_{R,N}(t,\Psi^{-1}(y)) && (t,y)\in(-T_0N^{-2},T_0N^{-2})\times(\mathbb{R}\times\mathbb{T}^3). \end{aligned}$$

*Then there is* $R_0\geq 1$ *(depending on* $T_0$, $\phi'$ *and* $\varepsilon_1$*) such that, for any* $R\geq R_0$ *and* $N\geq 10R$,

$$\limsup_{N\to\infty}\|U_N-V_{R,N}\|_{X^1(-T_0N^{-2},T_0N^{-2})}\lesssim_{E_{\mathbb{R}^4}(\phi)}\varepsilon_1.$$

*Proof.* We show **Part (i)** and **Part (ii)** together, by Proposition 4.5. By using Theorem 5.1, we know that $v'$ exists globally and satisfies

$$\|\nabla_{\mathbb{R}^4}v'\|_{(L_t^\infty L_x^2\cap L_t^2L_x^4)(\mathbb{R}\times\mathbb{R}^4)}\lesssim 1,$$

and

$$\sup_{t\in\mathbb{R}}\|v'(t)\|_{H^5(\mathbb{R}^4)}\lesssim_{\|\phi'\|_{H^5(\mathbb{R}^4)}}1. \tag{5.4}$$

Recall that $v'_R(t,x)=\eta(x/R)v'(t,x)$, which satisfies the following equation:

$$\begin{aligned} &(i\partial_t+\Delta_{\mathbb{R}^4})v'_R=(i\partial_t+\Delta_{\mathbb{R}^4})(\eta(x/R)v'(t,x))\\ &=\eta(x/R)(i\partial_t+\Delta_{\mathbb{R}^4})v'(t,x)+R^{-2}v'(t,x)(\Delta_{\mathbb{R}^4}\eta)(x/R)+2R^{-1}\sum_{j=1}^4\partial_jv'(t,x)\partial_j\eta(x/R). \end{aligned}$$

Then it implies

$$(i\partial_t+\Delta_{\mathbb{R}^4})v'_R=\rho\,|v'_R|^2v'_R+e_R(t,x),$$

where

$$e_R(t,x)=\rho(\eta(x/R)-\eta^3(x/R))v'|v'|^2+R^{-2}v'(t,x)(\Delta_{\mathbb{R}^4}\eta)(x/R)+2R^{-1}\sum_{j=1}^4\partial_jv'(t,x)\partial_j\eta(x/R).$$

After scaling, we have that

$$(i\partial_t+\Delta_{\mathbb{R}^4})v'_{R,N}=\rho\,|v'_{R,N}|^2v'_R+e_{R,N}(t,x),$$

where $e_{R,N}(t,x)=N^3e_R(N^2t,Nx)$ with $V_{R,N}(t,y)=v'_{R,N}(t,\Psi^{-1}(y))$ and taking $N\geq 10R$, we obtain that

$$(i\partial_t+\Delta_{\mathbb{R}^4})V_{R,N}(t,y)=\rho\,|V_{R,N}|^2V_{R,N}+E_{R,N}(t,y), \tag{5.5}$$

where $E_{R,N}(t,y)=e_{R,N}(t,\Psi^{-1}(y))$. In order to apply Theorem 4.5, we need to verify the following three conditions:

**(1)** $\|V_{R,N}\|_{L_t^\infty([-T_0N^{-2},T_0N^{-2}]:H^1(\mathbb{R}\times\mathbb{T}^3))}+\|V_{R,N}\|_{Z([-T_0N^{-2},T_0N^{-2}])}\leq M$;
**(2)** $\|f_N-V_{R,N}(0)\|_{H^1(\mathbb{R}\times\mathbb{T}^3)}\leq\varepsilon$;
**(3)** $\|E_{R,N}\|_{N([-T_0N^{-2},T_0N^{-2}])}\leq\varepsilon$.

We now verify these conditions one by one.

**Condition 1:** $\|V_{R,N}\|_{L_t^\infty([-T_0N^{-2},T_0N^{-2}]:H^1(\mathbb{R}\times\mathbb{T}^3))}+\|V_{R,N}\|_{Z([-T_0N^{-2},T_0N^{-2}])}\leq M$.
Since $v'(t,x)$ globally exists, $V_{R,N}(t,y)$ also globally exists. For any given $T_0\in(0,\infty)$, we have

$$\begin{aligned} &\sup_{t\in[-T_0N^{-2},T_0N^{-2}]}\|V_{R,N}(t)\|_{H^1(\mathbb{R}\times\mathbb{T}^3)}\leq\sup_{t\in[-T_0N^{-2},T_0N^{-2}]}\|v'_{R,N}(t)\|_{H^1(\mathbb{R}^4)}\\ &=\sup_{t\in[-T_0N^{-2},T_0N^{-2}]}\frac{1}{N}\|v'_R(N^2t)\|_{L^2(\mathbb{R}^4)}+\|v'_R(N^2t)\|_{\dot H^1(\mathbb{R}^4)}\\ &\lesssim\|v'(t,x)\|_{H^1(\mathbb{R}^4)}\leq\|\phi'(t)\|_{H^5(\mathbb{R}^4)}. \end{aligned}$$

By Littlewood-Paley theorem and Sobolev embedding, we obtain that

$$\begin{aligned}
&\|V_{R,N}\|_{Z([-T_0N^{-2},T_0N^{-2}])}\\
\leq &\sup_{J\subset[-T_0N^{-2},T_0N^{-2}]}\|(\sum_M(\langle\nabla\rangle^{\frac{1}{2}}P_MV_{R,N})^2)^{\frac{1}{2}}\|_{L^4_{t,x}(J\times\mathbb{R}\times\mathbb{T}^3)}\\
\lesssim &\sup_{J\subset[-T_0N^{-2},T_0N^{-2}]}\|\langle\nabla\rangle^{\frac{1}{2}}V_{R,N}\|_{L^4_{t,x}(J\times\mathbb{R}\times\mathbb{T}^3)}\\
\leq &\sup_{J\subset[-T_0N^{-2},T_0N^{-2}]}\|\langle\nabla\rangle V_{R,N}\|_{L^4_tL^{\frac{8}{3}}_x(J\times\mathbb{R}\times\mathbb{T}^3)}\\
\lesssim &\|v'_R\|_{L^4_tL^{\frac{8}{3}}_x([-T_0,T_0]\times\mathbb{R}^4)}+\|\nabla_{\mathbb{R}^4}v'_R\|_{L^4_tL^{\frac{8}{3}}_x([-T_0,T_0]\times\mathbb{R}^4)}.
\end{aligned}$$

Since

$$\|v'_R\|_{L^4_tL^{\frac{8}{3}}_x([-T_0,T_0]\times\mathbb{R}^4)}+\|\nabla_{\mathbb{R}^4}v'_R\|_{L^4_tL^{\frac{8}{3}}_x([-T_0,T_0]\times\mathbb{R}^4)}\lesssim\sup_t\|v'(t)\|_{H^5}$$

and by (5.4), we obtain

$$\|V_{R,N}\|_{Z([-T_0N^{-2},T_0N^{-2}])}\lesssim_{\|\phi'\|_{H^5(\mathbb{R}^4)}}1.$$

**Condition 2:** $\|f_N-V_{R,N}(0)\|_{H^1(\mathbb{R}\times\mathbb{T}^3)}\leq\varepsilon$.
By Hölder inequality, we have that

$$\begin{aligned}
\|f_N-V_{R,N}(0)\|_{H^1(\mathbb{R}\times\mathbb{T}^3)}&\leq\|\phi_N(\Psi^{-1}(y))-\phi'_{R,N}(\Psi^{-1}(y))\|_{\dot{H}^1(\mathbb{R}\times\mathbb{T}^3)}\\
&\leq\|\eta(x/N^{\frac{1}{2}})\phi(x)-\phi(x)\|_{\dot{H}^1(\mathbb{R}^4)}+\|\phi-\phi'\|_{\dot{H}^1(\mathbb{R}^4)}+\|\eta(x/N^{\frac{1}{2}})\phi'(x)-\phi'(x)\|_{\dot{H}^1(\mathbb{R}^4)}.
\end{aligned}$$

With $N\geq 10R$, and $R>R_0$, $R_0$ large enough, we obtain that

$$\|f_N-V_{R,N}(0)\|_{H^1(\mathbb{R}\times\mathbb{T}^3)}\leq 2\varepsilon_1.$$

**Condition 3:** $\|E_{R,N}\|_{N([-T_0N^{-2},T_0N^{-2}])}\leq\varepsilon$.
By Proposition 2.8 and the scaling invariance, we write that

$$\begin{aligned}
\|E_{R,N}\|_{N([-T_0N^{-2},T_0N^{-2}])}&=\left\|\int_0^t e^{i(t-s)\Delta}E_{R,N}(s)\,ds\right\|_{X^1([-T_0N^{-2},T_0N^{-2}])}\\
&\leq\sup_{\|u_0\|_{Y^{-1}}=1}\|u_0\|_{Y^{-1}}\||\nabla|E_{R,N}\|_{L^1_tL^2_x([-T_0N^{-2},T_0N^{-2}]\times\mathbb{R}\times\mathbb{T}^3)}\\
&=\|\nabla_{\mathbb{R}^4}e_R\|_{L^1_tL^2_x([-T_0,T_0]\times\mathbb{R}^4)}.
\end{aligned}$$

For the $|\nabla_{\mathbb{R}^4}e_R(t,x)|$, we have the following estimate:

$$\begin{aligned}
&|\nabla_{\mathbb{R}^4}e_R(t,x)|\\
&\leq\left|\nabla_{\mathbb{R}^4}(\eta(x/R)-\eta(x/R)^3)v'(t,x)|v'(t,x)|^2\right|+3\left|(\eta(x/R)-\eta(x/R)^3)\nabla_{\mathbb{R}^4}v'(t,x)|v'(t,x)|^2\right|\\
&\quad+R^{-3}\left|v'(t,x)\nabla_{\mathbb{R}^4}\Delta_{\mathbb{R}^4}\eta(x/R)\right|+R^{-2}\left|\nabla_{\mathbb{R}^4}v'(t,x)(\Delta_{\mathbb{R}^4}\eta)(x/R)\right|+R^{-1}\Big|\sum_{i,j=1}^4\partial_{i,j}v'\nabla_{\mathbb{R}^4}\eta(x/R)\Big|\\
&\lesssim_{\|\phi'\|_{H^5(\mathbb{R}^4)}}\mathbb{1}_{[R,2R]}(|x|)\left(|v'(t,x)|+|\nabla_{\mathbb{R}^4}v'(t,x)|\right)+\frac{1}{R}\left(|\langle\nabla_{\mathbb{R}^4}\rangle^2v'(t,x)|\right).
\end{aligned}$$

Since $\|\nabla^2_{\mathbb{R}^4}v'\|_{L^\infty_x}\lesssim_{\|\phi'\|_{H^5(\mathbb{R}^4)}}1$, $\|\sum_{i,j=1}^4\partial_{i,j}v'\|_{L^\infty_x}\lesssim_{\|\phi'\|_{H^5(\mathbb{R}^4)}}1$ and $\|v'\|_{L^\infty_x}\lesssim_{\|\phi'\|_{H^5(\mathbb{R}^4)}}1$ (by Sobolev embedding), we obtain that

$$\begin{aligned}
&\|\nabla_{\mathbb{R}^4}e_R\|_{L^1_tL^2_x([-T_0,T_0]\times\mathbb{R}^4)}=\int_{-T_0}^{T_0}(\int_{\mathbb{R}^4}|\nabla_{\mathbb{R}^4}e_R|^2\,dx)^{\frac{1}{2}}\,dt\\
&\leq\int_{-T_0}^{T_0}\left(\int_{\mathbb{R}^4}\mathbb{1}_{[R,2R]}(|x|)(|v'(t,x)|^2+|\nabla_{\mathbb{R}^4}v'(t,x)|^2)\,dx+\frac{1}{R^2}\int_{\mathbb{R}^4}|\langle\nabla_{\mathbb{R}^4}\rangle^2v'(t,x)|^2\,dx\right)^{\frac{1}{2}}dt\\
&\lesssim_{\|\phi'\|_{H^5}}2T_0\left(\int_{\mathbb{R}^4}\mathbb{1}_{[R,2R]}(|x|)\langle\nabla_{\mathbb{R}^4}\rangle^2v'(t,x)|^2\,dx\right)^{\frac{1}{2}}+\frac{1}{R}\to 0,\ \text{as } R\to\infty.
\end{aligned}$$

So we can obtain that

$$\|\nabla E_{R,N}\|_{L_t^1L_x^2([-T_0N^{-2},T_0N^{-2}]\times\mathbb{R}\times\mathbb{T}^3)}<\varepsilon_1,$$

where $R>R_0$, and $R_0$ large enough.

All the three conditions checked, we have the desired result. □

Next, we will introduce another useful extinction lemma, which is proved in [22].

**Lemma 5.3** (Extinction Lemma, Lemma 4.3 in [22])**.** *Suppose that $\phi\in\dot H^1(\mathbb{R}^4),\varepsilon>0$, and an interval $I\subset\mathbb{R}$. Assume that*

$$\|\phi\|_{\dot H^1(\mathbb{R}^4)}\le 1,\quad \|\nabla e^{it\Delta}\phi\|_{L_t^2L_x^4(I\times\mathbb{R}^4)}\le\varepsilon.$$

*For $N\ge 1$, we define as before (recall the definitions in* (5.3)*):*

$$Q_N\phi=\eta(N^{-1/2}x)\phi(x),\quad \phi_N=N(Q_N\phi)(Nx),\quad f_N(y)=\phi_N(\Psi^{-1}(y)).$$

*Then there exists $N_0=N_0(\phi,\varepsilon)$ such that for any $N\ge N_0$,*

$$\|e^{it\Delta}f_N\|_{Z(N^{-2}I)}\lesssim\varepsilon.$$

**Definition 5.4.** Let $\widetilde{\mathcal{F}}_e$ denote the set of *renormalized Euclidean frames*

$$\begin{aligned}\widetilde{\mathcal{F}}_e:=&\{(N_k,t_k,x_k)_{k\ge1}:N_k\in[1,\infty),\ t_k\in[-1,1],\ x_k\in\mathbb{R}\times\mathbb{T}^3,\ N_k\to\infty\\ &\text{and either } t_k=0 \text{ for any } k\ge1 \text{ or } \lim_{k\to\infty}N_k^2|t_k|=\infty\}.\end{aligned}$$

Let us introduce a useful proposition that shows that the solutions based on Euclidean profiles evolve as the solutions in the the Euclidean space which scatter.

**Proposition 5.5.** *Assume that $\mathcal{O}=(N_k,t_k,x_k)_k\in\widetilde{\mathcal{F}}_e$, $\phi\in\dot H^1(\mathbb{R}^4)$, and let $U_k(0)=\Pi_{t_k,x_k}(T_{N_k}\phi)$. Under the assumption that for any solutions $v$ of* (5.1) *with initial datum $\phi$, we have*

$$\sup_{t\in \text{lifespan of } v}\|v(t)\|_{\dot H^1(\mathbb{R}^4)}<\|W\|_{\dot H^1(\mathbb{R}^4)}.$$

*Then the following two conclusions hold:*

**Part (i)** *For $k$ large enough (depending only on $\phi$ and $\mathcal{O}$), there is a nonlinear solution to* (1.1) *$U_k\in X^1([-1,1])$ with initial datum $U_k(0)$ satisfying:*

$$\|U_k\|_{X^1([-1,1])}\lesssim_{E_{\mathbb{R}^4}(\phi)}1. \tag{5.6}$$

**Part (ii)** *There exists a Euclidean solution $u\in C(\mathbb{R}:\dot H^1(\mathbb{R}^4))$ of*

$$(i\partial_t+\Delta_{\mathbb{R}^4})u=-|u|^2u$$

*with scattering data $\phi^{\pm\infty}$ defined as in* (5.2) *such that up to a subsequence: for any $\varepsilon>0$, there exists $T(\phi,\varepsilon)$ such that for all $T\ge T(\phi,\varepsilon)$ there exists $R(\phi,\varepsilon,T)$ such that for all $R\ge R(\phi,\varepsilon,T)$, there holds that*

$$\|U_k-\widetilde u_k\|_{X^1(\{|t-t_k|\le TN_k^{-2}\})}\le\varepsilon, \tag{5.7}$$

*for $k$ large enough, where*

$$(\pi_{-x_k}\widetilde u)(t,x)=N_k\eta(N_k\Psi^{-1}(x)/R)u(N_k^2(t-t_k),N_k\Psi^{-1}(x)).$$

*In addition, up to a subsequence,*

$$\|U_k(t)-\Pi_{t_k-t,x_k}T_{N_k}\phi^{\pm\infty}\|_{X^1(\{\pm(t-t_k)\ge TN_k^{-2}\})}\le\varepsilon$$

*for $k$ large enough (depending on $\phi,\varepsilon,T,R$). Recall definitions of $\pi$ and $\Pi$ in* (4.1).

*Proof.* Without loss of generality, it is sufficient to prove the case $x_k = 0$.

First we consider **Part (i)**. For $k$ large enough, we can make

$$\|\phi - \eta(x/N_k^{\frac{1}{2}})\phi\|_{\dot{H}^1(\mathbb{R}^4)} \leq \varepsilon_1$$

for any $\varepsilon_1 > 0$. For each $N_k$, we choose $T_{0,N_k} = N_k^2$ (Recall $T_{0,N_k}$ in Lemma 5.2). For each $T_{0,N_k}$, we can choose $R_k$ large enough to make it possible to apply Lemma 5.2. Note that in this case, the choice of $R_k$ is determined by $T_{0,N_k}$ as in the proof of Lemma 5.2, as a result Lemma 5.2 gives $U_k$ satisfying (5.6).

We then consider **Part (ii)**. The first case in Euclidean frame is when $t_k = 0$ for all $k$. To prove (5.7), we need to choose $T(\phi, \delta)$ large enough, to make sure

$$\|\nabla_{\mathbb{R}^4} u\|_{L^3_{t,x}(\{|t|>T(\phi,\delta)\}\times\mathbb{R}^4)} \leq \delta.$$

By Theorem 5.1, we obtain that

$$\|u(\pm T(\phi,\delta)) - e^{\pm iT(\phi,\delta)\Delta}\phi^{\pm\infty}\|_{\dot{H}^1(\mathbb{R}^4)} \leq \delta,$$

which implies

$$\|U_{N_k}(\pm TN_k^{-2}) - \Pi_{\mp T,x_k} T_{N_k}\phi^{\pm\infty}\|_{H^1(\mathbb{R}\times\mathbb{T}^3)} \leq \delta.$$

Based on the definition of $Z$-norm and Proposition 4.1, we have

$$\|e^{it\Delta}\left(U_{N_k}(\pm TN_k^{-2}) - \Pi_{\mp T,x_k} T_{N_k}\phi^{\pm\infty}\right)\|_{X^1(|t|<T^{-1})} \lesssim \delta. \tag{5.8}$$

Then by Theorem 4.3, we obtain that

$$\|U_{N_k} - e^{it\Delta}U_{N_k}(\pm TN_k^{-2})\|_{X^1} \leq \delta, \tag{5.9}$$

and combining (5.8) and (5.9), we have

$$\|U_{N_k} - \Pi_{-t,x_k} T_{N_k}\phi^{\pm\infty}\|_{X^1(\{\pm t \geq \pm TN_k^{-2}\}\cap\{|t|<T^{-1}\})} \leq \varepsilon,$$

where we choose $\delta$ small enough.

The second case is when $N_k^2|t_k| \to \infty$. In this case, we have

$$U_k(0) = \Pi_{t_k,0}(T_{N_k}\phi) = e^{-it_k\Delta}\left(N_k^{\frac{1}{2}}\eta(N_k^{\frac{1}{2}}\Psi^{-1}(x))\phi(N_k\Psi^{-1}(x))\right).$$

Because of the existence of wave operators of NLS based on Theorem 5.1, we know the following final state value problem

$$\begin{cases} (i\partial_t + \Delta_{\mathbb{R}^4})v = -v|v|^2, \\ \lim_{t\to-\infty} \|v(t) - e^{it\Delta}\phi\|_{\dot{H}^1(\mathbb{R}^4)} = 0, \end{cases}$$

is global well-posed under the assumption of the *a priori* bound $\sup_{t\in\text{lifespan of } v} \|v(t)\|_{\dot{H}^1(\mathbb{R}^4)} < \|W\|_{\dot{H}^1(\mathbb{R}^4)}$. We then set that

$$\widetilde{v_k}(t) := N_k^{\frac{1}{2}}\eta(N_k\Psi^{-1}(x)/R)v(N_k\Psi^{-1}(x), N_k^2 t),$$

so we have $\widetilde{v_k}(-t_k) = N_k^{\frac{1}{2}}\eta(N_k\Psi^{-1}(x)/R)v(N_k\Psi^{-1}(x), -N_k^2 t_k)$. For $R$ large enough, we obtain that

$$\|\widetilde{v_k}(-t_k) - U_k(0)\|_{H^1(\mathbb{R}\times\mathbb{T}^3)} \leq \|\eta(x/N_k^{\frac{1}{2}})v(x, -N_k^2 t_k) - e^{it_k N_k^2\Delta}\eta(x/N_k^{\frac{1}{2}})\phi(x)\|_{H^1(\mathbb{R}^4)} \to 0,$$

as $k \to \infty$. Denote that $V_k(t) := \widetilde{v_k}(t-t_k)$ and then as above we have $\lim_{k\to\infty}\|U_k(0)-V_k(0)\|_{H^1(\mathbb{R}\times\mathbb{T}^3)} = 0$. Hence by the stability property (Theorem 4.5), we have that for $k$ large enough $U_k$ exists and satisfies

$$\|U_k - V_k\|_{X^1([-1,1])} \to 0, \text{ as } k \to \infty. \tag{5.10}$$

By the definition of $V_k$ and (5.10), we then can prove all desired properties of $U_k$. $\square$

**Lemma 5.6.** *Assume that* $\mathcal{O}=(N_k,t_k,x_k)_k\in\widetilde{\mathcal{F}}_e$, $\phi\in\dot{H}^1(\mathbb{R}^4)$. *Let* $U_k(0)=\Pi_{t_k,x_k}(T_{N_k}\phi)$ *and* $U_k$ *is the solution of* (1.1) *with the lifespan* $I_k=[-T_k,T^k]$. *Assume also that there exists* $C>0$ *such that*

$$\limsup_{k\to\infty}\sup_{t\in[-T_k,T^k]}\|U_k(t)\|_{H^1(\mathbb{R}\times\mathbb{T}^3)}<C<\|W\|_{\dot{H}^1(\mathbb{R}^4)}. \tag{5.11}$$

*Denote that* $T_\infty=\limsup_{k\to\infty}N_k^2T_k$ *and* $T^\infty=\limsup_{k\to\infty}N_k^2T^k$. *Then if* $T_\infty=T^\infty=\infty$, *under the assumption that for any solution* $v$ *of* (5.1) *with initial datum* $\phi$, *we have*

$$\sup_{t\in \text{lifespan of } v}\|v(t)\|_{\dot{H}^1(\mathbb{R}^4)}<\|W\|_{\dot{H}^1(\mathbb{R}^4)}. \tag{5.12}$$

*Proof of Lemma 5.6.* Denote that $\widetilde{v}_k:=\Pi_{t_k,x_k}(T_{N_k}v)$. By **Part (ii)** in Proposition 5.2 (taking $N=N_k$ and $R=N_k^{\frac{1}{2}}$) and $N_k\to\infty$ as $k\to\infty$, we have

$$\limsup_{k\to\infty}\sup_{t\in[-T_k,T^k]}\|U_k(t)-\widetilde{v}_k(t)\|_{H^1(\mathbb{R}\times\mathbb{T}^3)}\lesssim\limsup_{k\to\infty}\|U_k-\widetilde{v}_k\|_{X^1([-T_k,T^k])}=0,$$

which together with (5.11) implies

$$\limsup_{k\to\infty}\sup_{t\in[-T_k,T^k]}\|\widetilde{v}_k(t)\|_{H^1(\mathbb{R}\times\mathbb{T}^3)}<\|W\|_{\dot{H}^1(\mathbb{R}^4)}.$$

By taking $N_k\to\infty$, we then have (5.12). □

**Proposition 5.7** (Decomposition of the nonlinear Euclidean profiles). *Assume that* $U_k$ *is the nonlinear Euclidean profile associated to* $\mathcal{O}=(N_k,t_k,x_k)_k\in\widetilde{\mathcal{F}}_e$ *defined as in Lemma 5.6 and* $U_k$ *also satisfies the condition* (5.11). *For any* $\theta>0$, *there exists a large enough* $T_0^\theta$, *such that for any* $T^\theta\geq T_0^\theta$, $R^\theta$ *sufficiently large and all* $k$ *large enough (depending on* $R^\theta$*) we can decompose* $U_k$ *into*

$$U_k=\omega_k^{\theta,-\infty}+\omega_k^{\theta,+\infty}+\omega_k^\theta+\rho_k^\theta. \tag{5.13}$$

*Here the components* $\omega_k^{\theta,\pm\infty}$, $\omega_k^\theta$, *and* $\rho_k^\theta$ *satisfy the following conditions:*

$$\begin{aligned}
&\|\omega_k^{\theta,\pm\infty}\|_{Z([-1,1])}+\|\rho_k^\theta\|_{X^1([-1,1])}+\|(1-P_{\leq R^\theta N_k}+P_{\leq(R^\theta)^{-1}N_k})\,\omega_k^\theta\|_{X^1([-1,1])}\leq\theta,\\
&\|\omega_k^{\theta,\pm\infty}\|_{X^1([-1,1])}+\|\omega_k^\theta\|_{X^1([-1,1])}\lesssim1,\\
&\omega_k^{\theta,\pm\infty}=\mathbb{1}_{\pm(t-t_k)\geq T^\theta N_k^{-2}}(t)\,\Pi_{(t_k-t,x_k)}T_{N_k}\phi^{\pm\infty,\theta},\quad \|\phi^{\pm\infty,\theta}\|_{\dot{H}^1}\lesssim1,\\
&|\nabla_x^m\omega_k^\theta|+(N_k)^{-2}\mathbb{1}_{S_k^\theta}|\partial_t\nabla_x^m\omega_k^\theta|\leq R^\theta(N_k)^{|m|+1}\mathbb{1}_{S_k^\theta},\quad m\in\mathbb{N}^4,\ 0\leq|m|\leq10,
\end{aligned} \tag{5.14}$$

*where*

$$S_k^\theta:=\{(t,x)\in[-1,1]\times\mathbb{R}\times\mathbb{T}^3:|t-t_k|<T^\theta(N_k)^{-2},\ |x-x_k|\leq R^\theta(N_k)^{-1}\}.$$

*Proof of Proposition 5.7.* Recall that we defined $T_\infty=\limsup_{k\to\infty}N_k^2T_k$ and $T^\infty=\limsup_{k\to\infty}N_k^2T^k$ in Lemma 5.6. Now we will consider the following two cases separately.

**(1)** $T_\infty=T^\infty=\infty$,
**(2)** At least one of $T_\infty$ and $T^\infty$ is finite.

**Case 1:** First we consider $T_\infty=T^\infty=\infty$.
In this case, we can construct the decomposition (5.13) using Lemma 5.6 and Proposition 5.5. By Proposition 5.5, there exists $T(\phi,\frac{\theta}{4})$, such that for all $T\geq T(\phi,\frac{\theta}{4})$, there exists $R(\phi,\frac{\theta}{4},T)$ such that for all $R\geq R(\phi,\frac{\theta}{2},T)$, there holds that

$$\|U_k-\widetilde{u}_k\|_{X^1(\{|t-t_k|\leq T(N_k)^{-2}\}\cap\{|t|\leq1\})}\leq\frac{\theta}{2},$$

for $k$ large enough, where

$$(\pi_{-x_k}\widetilde{u}_k)\,(t,x)=N_k\eta(N_k\Psi^{-1}(x)/R)u(N_k^2(t-t_k),N_k\Psi^{-1}(x)), \tag{5.15}$$

for $u$ being the solution of (5.1) with scattering data $\phi^{\pm\infty}$.

In addition, up to subsequence, we have

$$\|U_k - \Pi_{t_k - t, x_k} T_{N_k} \phi^{\pm\infty}\|_{X^1(\{\pm(t-t_k)\geq T(N_k)^{-2}\}\cap\{|t|\leq 1\})} \leq \frac{\theta}{4},$$

for $k$ large enough (depending on $\phi$, $\theta$, $T$, and $R$). Then we choose a sufficiently large $T^\theta > T(\phi, \frac{\theta}{4})$ based on the extinction lemma (Lemma 5.3), such that

$$\|e^{it\Delta}\Pi_{t_k,x_k} T_{N_k}\phi^{\pm\infty}\|_{Z(T^\theta(N_k)^{-2},(T^\theta)^{-1})} \leq \frac{\theta}{4}$$

for $k$ large enough. Accordingly, we define $R^\theta = R(\phi, \frac{\theta}{2}, T^\theta)$.

Now we are ready to give the decomposition in (5.13). We take the following choices as the components in (5.13) and we will also verify that these choices actually satisfy all the requirements in (5.14).

**(1)** $\omega_k^{\theta,\pm\infty} := \mathbb{1}_{\{\pm(t-t_k)\geq T^\theta(N_k)^{-2},|t|\leq 1\}} \left(\Pi_{t_k-t,x_k} T_{N_k}\phi^{\theta,\pm\infty}\right)$
Here $\phi^{\theta,\pm\infty} = P_{\leq R^\theta}(\phi^{\theta,\pm\infty})$ and $\|\phi^{\theta,\pm\infty}\|_{\dot{H}^1(\mathbb{R}^4)} \lesssim 1$. This choice implies $\omega_k^{\theta,\pm\infty} = P_{\leq R^\theta N_k}\omega_k^{\theta,\pm\infty}$.

**(2)** $\omega_k^\theta := \widetilde{u}_k \cdot \mathbb{1}_{S_k^\theta}$
By the stability property (Theorem 4.5) and Theorem 5.2, we can adjust $\omega_k^\theta$ and $\omega_k^{\theta,\pm\infty}$ with an acceptable error, and make

$$|\nabla_x^m \omega_k^\theta| + (N_k)^{-2}\mathbb{1}_{S_k^{\alpha,\theta}}|\partial_t \nabla_x^m \omega_k^\theta| \leq R^\theta (N_k)^{|m|+1}\mathbb{1}_{S_k^\theta},\ 0 \leq |m| \leq 10.$$

**(3)** $\rho_k^\theta := U_k^\alpha - \omega_k^\theta - \omega^{\theta,+\infty} - \omega^{\theta,-\infty}$.
By the definitions of $\omega_k^{\theta,\pm\infty}$, $\omega_k^\theta$ and $\rho_k^\theta$ and Proposition 5.5, we obtain that

$$\|\rho_k^\theta\|_{X^1([-1,1])} \leq \frac{\theta}{2}.$$

and then we have

$$\|\omega_k^{\theta,\pm\infty}\|_{Z'([-1,1])} + \|\rho_k^\theta\|_{X^1([-1,1])} + \|(1 - P_{\leq R^\theta N_k} + P_{\leq (R^\theta)^{-1}N_k})\,\omega_k^\theta\|_{X^1([-1,1])} \leq \theta,$$
$$\|\omega_k^{\theta,\pm\infty}\|_{X^1([-1,1])} + \|\omega_k^\theta\|_{X^1([-1,1])} \lesssim 1.$$

**Case 2:** In the end, we treat $T_\infty < \infty$ or/and $T^\infty < \infty$.
It is easy to check that we can construct the decomposition in the same way as in **Case 1** but (5.13) by taking $\omega_k^{\theta,-\infty} = 0$ or/and $\omega_k^{\theta,-\infty} = 0$ respectively. □

**Corollary 5.8** (Decomposition of the nonlinear scale-one profiles)**.** *Assume that $U_k$ is the nonlinear scale-one profile associated to $\mathcal{O} = (1, t_k, x_k)_k \in \mathcal{F}_1$ with initial datum $U_k(0) = \widetilde{\psi}_{\mathcal{O}_k^\alpha}$ and $U_k$ satisfies the condition* (5.11)*. For any $\theta > 0$, there exists a large enough $T_0^\theta$, such that for any $T^\theta \geq T_0^\theta$, $R^\theta$ sufficiently large and all $k$ large enough (depending on $R^\theta$), we can decompose $U_k$ as what we did in Proposition 5.7*

$$U_k = \omega_k^{\theta,-\infty} + \omega_k^{\theta,+\infty} + \omega_k^\theta + \rho_k^\theta,$$

*where $\rho_k^\theta = \omega_k^\theta = 0$ and $\omega^{\theta,+\infty} = \omega^{\theta,-\infty}$ and the components $\omega_k^{\theta,\pm\infty}$, $\omega_k^\theta$, and $\rho_k^\theta$ satisfy the same requirements as in* (5.14)*.*

*Proof of Corollary 5.8.* By taking $T_0^\theta$ large, it is easy to make $\|\omega^{\theta,\pm\infty}\|_{Z'(-T_0^\theta,T_0^\theta)} \leq \theta$. Hence the decomposition for both types of nonlinear profiles (Euclidean and scale-one) can be unified. □

**Lemma 5.9** (Almost orthogonality of nonlinear profiles)**.** *Define $U_k^\alpha$, $U_k^\beta$ to be the solutions to* (1.1) *with the life-span $I_k$ and initial data $U_k^\alpha(0) = \widetilde{\psi}^\alpha_{\mathcal{O}_k^\alpha}$, $U_k^\beta(0) = \widetilde{\psi}^\beta_{\mathcal{O}_k^\beta}$, where $\psi^\alpha, \psi^\beta \in H^1(\mathbb{R}\times\mathbb{T}^3)$ and frames $\mathcal{O}^\alpha$ and $\mathcal{O}^\beta$ are orthogonal. Then*

$$\lim_{k\to\infty} \sup_{t\in I_k} \langle U_k^\alpha(t), U_k^\beta(t)\rangle_{\dot{H}^1\times\dot{H}^1} = 0.$$

*Proof of Lemma 5.9.* By the unified decomposition in Proposition 5.7 and Corollary 5.8, we know that for the nonlinear profiles $U^\alpha$ and $U^\beta$, there exist $T^{\alpha,\theta}$, $R^{\alpha,\theta}$, $T^{\beta,\theta}$, $R^{\beta,\theta}$ sufficiently large such that

$$\begin{aligned}U_k^\alpha &= \omega_k^{\alpha,\theta,-\infty}+\omega_k^{\alpha,\theta,+\infty}+\omega_k^{\alpha,\theta}+\rho_k^{\alpha,\theta},\\ U_k^\beta &= \omega_k^{\beta,\theta,-\infty}+\omega_k^{\beta,\theta,+\infty}+\omega_k^{\beta,\theta}+\rho_k^{\beta,\theta}.\end{aligned}$$

Note that $\rho_k^{\alpha,\theta}$, $\rho_k^{\beta,\theta}$ are the small terms with their $X^1$ norm less then $\theta$, for any fixed $t\in I_k$. Then to show the almost orthogonality, it will suffice to consider the following three interactions:

**(1)** $\langle\omega_k^{\alpha,\theta,\pm\infty},\omega_k^{\beta,\theta,\pm\infty}\rangle_{\dot H^1\times\dot H^1}$;
**(2)** $\langle\omega_k^{\alpha,\theta,\pm\infty},\omega_k^{\beta,\theta}\rangle_{\dot H^1\times\dot H^1}$;
**(3)** $\langle\omega_k^{\alpha,\theta},\omega_k^{\beta,\theta}\rangle_{\dot H^1\times\dot H^1}$.

**Case 1:** $\langle\omega_k^{\alpha,\theta,\pm\infty},\omega_k^{\beta,\theta,\pm\infty}\rangle_{\dot H^1\times\dot H^1}$.
By the constructions of $\omega_k^{\alpha,\theta,\pm\infty},\omega_k^{\beta,\theta,\pm\infty}$ in the proof of Proposition 5.7, we write

$$\begin{aligned}\omega_k^{\alpha,\theta,\pm\infty} &:= \mathbb{1}_{\{\pm(t-t_k^\alpha)\geq T^{\alpha,\theta}(N_k^\alpha)^{-2}\}}\left(\Pi_{t_k^\alpha-t,x_k^\alpha}T_{N_k^\alpha}\phi^{\alpha,\theta,\pm\infty}\right),\\ \omega_k^{\beta,\theta,\pm\infty} &:= \mathbb{1}_{\{\pm(t-t_k^\beta)\geq T^{\beta,\theta}(N_k^\beta)^{-2}\}}\left(\Pi_{t_k^\beta-t,x_k^\beta}T_{N_k^\beta}\phi^{\beta,\theta,\pm\infty}\right).\end{aligned}$$

For any fixed $t\in I_k$, we obtain that

$$\langle\omega_k^{\alpha,\theta,\pm\infty}(t),\omega_k^{\beta,\theta,\pm\infty}(t)\rangle_{\dot H^1\times\dot H^1}=\langle\phi_{\mathcal{O}_k^\alpha}^{\alpha,\theta,\pm\infty},\phi_{\mathcal{O}_k^\beta}^{\beta,\theta,\pm\infty}\rangle_{\dot H^1\times\dot H^1}.$$

We then arrive at

$$\lim_{k\to\infty}\sup_t\langle\omega_k^{\alpha,\theta,\pm\infty}(t),\omega_k^{\beta,\theta,\pm\infty}(t)\rangle_{\dot H^1\times\dot H^1}=0.$$

**Case 2:** $\langle\omega_k^{\alpha,\theta,\pm\infty},\omega_k^{\beta,\theta}\rangle_{\dot H^1\times\dot H^1}$.
By the constructions of $\omega_k^{\alpha,\theta,\pm\infty},\omega_k^{\beta,\theta,\pm\infty}$ again, we write

$$\omega_k^{\beta,\theta}:=\widetilde{u}_k^\beta\cdot\mathbb{1}_{S_k^{\beta,\theta}},$$

where $S_k^{\beta,\theta}:=\{(t,x)\in[-1,1]\times\mathbb{T}^4:|t-t_k^\beta|<T^{\beta,\theta}(N_k^\beta)^{-2},\ |x-x_k^\beta|\leq R^{\beta,\theta}(N_k^\beta)^{-1}\}$ and $\widetilde{u}_k^\beta$ is defined in (5.15). We then have the desired limit $\lim_{k\to\infty}\sup_t\langle\omega_k^{\alpha,\theta,\pm\infty},\omega_k^{\beta,\theta}\rangle_{\dot H^1\times\dot H^1}=0$.

**Case 3:** $\langle\omega_k^{\alpha,\theta},\omega_k^{\beta,\theta}\rangle_{\dot H^1\times\dot H^1}$.
We first take a small $\varepsilon>0$.

If $N_k^\alpha/N_k^\beta+N_k^\beta/N_k^\alpha\leq\varepsilon^{-1000}$ for $k$ large enough, then we claim that $S_k^{\alpha,\theta}\cap S_k^{\beta,\theta}=\emptyset$. In fact, by the definition of orthogonality of frames, $N_k^\alpha/N_k^\beta+N_k^\beta/N_k^\alpha\leq\varepsilon^{-1000}$ implies either $(N_k^\alpha)^2|t_k^\alpha-t_k^\beta|\to\infty$ or $N_k^\alpha|x_k^\alpha-x_k^\beta|\to\infty$, hence $S_k^{\alpha,\theta}\cap S_k^{\beta,\theta}=\emptyset$. Accordingly, the fact that there is not overlapping in their supports suggests $\omega_k^{\alpha,\theta}\,\omega_k^{\beta,\theta}\equiv 0$.

Otherwise, if $N_k^\alpha/N_k^\beta+N_k^\beta/N_k^\alpha>\varepsilon^{-1000}$, then without loss of generality, we assume that $N_k^\alpha/N_k^\beta\geq\varepsilon^{-1000}/2$. Using the same reasoning as in the previous scenario, we say that the supports of $\omega_k^{\alpha,\theta}$ and $\omega_k^{\beta,\theta}$ overlap and we can rewrite the product by rearranging the indicator function on their supports

$$\omega_k^{\alpha,\theta}\omega_k^{\beta,\theta}=\omega_k^{\alpha,\theta}\cdot(\omega_k^{\beta,\theta}\mathbb{1}_{(t_k^\alpha-T^{\alpha,\theta}(N_k^\alpha)^{-2},t_k^\alpha+T^{\alpha,\theta}(N_k^\alpha)^{-2})}(t))=:\omega_k^{\alpha,\theta}\widetilde{\omega}_k^{\beta,\theta}.$$

Now, in order to see the almost orthogonality, we decompose the interaction into several pieces and we claim that the interactions among all these terms are bounded by a multiple of $\varepsilon$, hence

$$\begin{aligned}\langle\omega_k^{\alpha,\theta},\omega_k^{\beta,\theta}\rangle_{\dot H^1\times\dot H^1}&\leq\langle P_{\leq\varepsilon^{10}N_k^\alpha}\omega_k^{\alpha,\theta},\widetilde{\omega}_k^{\beta,\theta}\rangle_{\dot H^1\times\dot H^1}+\langle P_{>\varepsilon^{10}N_k^\alpha}\omega_k^{\alpha,\theta},P_{>\varepsilon^{-10}N_k^\beta}\widetilde{\omega}_k^{\beta,\theta}\rangle_{\dot H^1\times\dot H^1}\\&\quad+\langle P_{>\varepsilon^{10}N_k^\alpha}\omega_k^{\alpha,\theta},P_{\leq\varepsilon^{-10}N_k^\beta}\widetilde{\omega}_k^{\beta,\theta}\rangle_{\dot H^1\times\dot H^1}\lesssim\varepsilon.\end{aligned}$$

In fact, the $\varepsilon$ bound above follows directly from $\varepsilon^{10}N_k^\alpha \gg \varepsilon^{-10}N_k^\beta$ and the following claim

**Claim 5.10.** For $k$ large enough and $\theta$ small enough, we have

**(1)** $\|\widetilde{\omega}_k^{\beta,\theta}\|_{X^1([-1,1])} \lesssim 1$;
**(2)** $\|P_{>\varepsilon^{-10}N_k^\beta}\widetilde{\omega}_k^{\beta,\theta}\|_{X^1([-1,1])} \lesssim \varepsilon$;
**(3)** $\|\omega_k^{\alpha,\theta}\|_{X^1([-1,1])} \lesssim 1$;
**(4)** $\|P_{\leq\varepsilon^{10}N_k^\alpha}\omega_k^{\alpha,\theta}\|_{X^1([-1,1])} \lesssim \varepsilon$.

Assuming Claim 5.10, we finish the proof of Lemma 5.9. Now we are left to prove Claim 5.10.

*Proof of Claim 5.10.* Denote $R := \max(R^{\alpha,\theta}, R^{\beta,\theta})$.

**(1)** First, let us consider $\|\widetilde{\omega}_k^{\beta,\theta}\|_{X^1([-1,1])}$. It is easy to check the following bound by Proposition 2.8

$$\begin{aligned}\|\widetilde{\omega}_k^{\beta,\theta}\|_{X^1([-1,1])} &\lesssim \left(\int_{|x-x_k^\beta|\leq R(N_k^\beta)^{-1}} |\langle\nabla\rangle\widetilde{\omega}_k^{\beta,\theta}(0)|^2\,dx\right)^{\frac{1}{2}} + \left(\sum_N \left(\int_{[-1,1]} \|P_N(\partial_t\widetilde{\omega}_k^{\beta,\theta})\|_{H^1} + \|P_N\Delta\widetilde{\omega}_k^{\beta,\theta}\|_{H^1}\,dt\right)^2\right)^{\frac{1}{2}}\\ &\lesssim R^3 + \int_{[-1,1]} (\|\partial_t\widetilde{\omega}_k^{\beta,\theta}\|_{H^1} + \|\Delta\widetilde{\omega}_k^{\beta,\theta}\|_{H^1})\,dt \quad \lesssim 1.\end{aligned}$$

**(2)** We then focus on the high frequency part of $\widetilde{\omega}_k^{\beta,\theta}$. Using Proposition 2.8 again, we have

$$\begin{aligned}&\|P_{>\varepsilon^{-10}N_k^\beta}\widetilde{\omega}_k^{\beta,\theta}\|_{X^1([-1,1])}\\ &\leq \left(\int_{|x-x_k^\beta|\leq R(N_k^\beta)^{-1}} |P_{>\varepsilon^{-10}N_k^\beta}\langle\nabla\rangle\widetilde{\omega}_k^{\beta,\theta}(0)|^2\,dx\right)^{\frac{1}{2}} + \int \|P_{>\varepsilon^{-10}N_k^\beta}(i\partial_t+\Delta)\widetilde{\omega}_k^{\beta,\theta}\|_{H^1}\,dt\\ &\leq \left(\int_{|x-x_k^\beta|\leq R(N_k^\beta)^{-1}} \left(\frac{\varepsilon^{10}}{N_k^\beta}\right)^2 |P_{>\varepsilon^{-10}N_k^\beta}\langle\nabla\rangle^2\widetilde{\omega}_k^{\beta,\theta}(0)|^2\,dx\right)^{\frac{1}{2}} + \int_{|t-t_k^\beta|<(N_k^\beta)^{-2}R} \frac{\varepsilon^{10}}{N_k^\beta}\|(i\partial_t+\Delta)\widetilde{\omega}_k^{\beta,\theta}\|_{H^2}\,dt\\ &\leq \varepsilon^{10}R^3 + (N_k^\beta)^{-2}R\frac{\varepsilon^{10}}{N_k^\beta}(R^4(N_k^\beta)^{-2}R^2(N_k^\beta)^{10})^{\frac{1}{2}} \quad \lesssim \varepsilon^{10}R^4.\end{aligned}$$

**(3)** The calculation for $\|\omega_k^{\alpha,\theta}\|_{X^1([-1,1])}$ is similar to what we did in **(1)**.
**(4)** Proposition 2.8 yields

$$\begin{aligned}\|P_{\leq\varepsilon^{10}N_k^\alpha}\omega_k^{\alpha,\theta}\|_{X^1([-1,1])} &\lesssim \varepsilon^{10}N_k^\alpha\left(\|P_{\leq\varepsilon^{10}N_k^\alpha}\omega_k^{\alpha,\theta}(0)\|_{L^2} + \int \|P_{\leq\varepsilon^{10}N_k^\alpha}(i\partial_t+\Delta)\omega_k^{\alpha,\theta}\|_{L^2}\,dt\right)\\ &\lesssim \varepsilon^{10}R^4.\end{aligned}$$

Now the proof of Claim 5.10 is complete. □

We finish the proof of Lemma 5.9. □

## 6. Rigidity Theorem

We are now ready to prove the main theorem by a contradiction argument. We follow the induction on energy method formalized in [24, 25]. Define the following functional

$$\Lambda(L) = \sup\{\|u\|_{Z(I)} : u \in X^1(I), \limsup_{t\in I}\|u(t)\|_{\dot{H}^1} \leq L, |I| \leq 1\}$$

where the supremum is taken over all strong solutions whose full energy is less than $L$. By the local theory, this functional is sublinear in $L$ and finite for $L$ sufficiently small. Note that this $\Lambda(L)$ is also non-decreasing. We define

$$L_{max} = \sup\{L : \Lambda(L) < +\infty\}. \tag{6.1}$$

Noticing the definition of $L_{max}$, showing Theorem 1.2 is equivalent to proving that $L_{max} \geq \|W\|_{\dot{H}^1}$ while in [22] their goal in the defocusing setting was to show $L_{max} = \infty$. This difference is originated from the different structure in the focusing dynamics and the possible existence of blow-up solutions beyond the ground state. The most important part in this global argument is to derive a contradiction assuming $L_{max} < \|W\|_{\dot{H}^1}$. It suffices to prove the key theorem as follows,

**Theorem 6.1.** *Consider $L_{max}$ defined in* (6.1)*, then $L_{max} \geq \|W\|_{\dot{H}^1(\mathbb{R}^4)}$.*

Clearly Theorem 6.1 implies the main Theorem 1.2. The main strategy to prove Theorem 6.1 is that, by using the profile decomposition, we analyze the following four scenarios respectively,

**(1)** no profiles
**(2)** exactly one Euclidean profile
**(3)** exactly one scale-one profile
**(4)** multiple (Euclidean/scale-one) profiles then rule out the possibility of all these cases to obtain the contradiction to the existence of such $L_{max}$. Notice that the first two cases are comparably easier. Now we prove Theorem 6.1 via profile decomposition.

*Proof of Theorem 6.1.* We prove $L_{max} \geq \|W\|_{\dot{H}^1}$ by a contradiction argument. First, we assume that $L_{max} < \|W\|_{\dot{H}^1}$. Then by the definition of $L_{max}$, there exists a sequence of solutions $u_k$ such that

$$\begin{aligned} &\limsup_{t\in[-T_k,T^k]} \|u_k(t)\|_{\dot{H}^1(\mathbb{R}\times\mathbb{T}^3)} \to L_{max}, \\ &\|u_k\|_{Z(-T_k,0)} \to +\infty, \\ &\|u_k\|_{Z(0,T^k)} \to +\infty, \end{aligned} \tag{6.2}$$

where $|T^k + T_k| \leq 1$.

Without loss of generality, we assume that $t_k = 0$. Then we define

$$L(\phi) = \limsup_{t\in[-T_k,T^k]} \|u_\phi(t)\|_{\dot{H}^1},$$

where $u_\phi(t)$ is the solution to (1.1) with initial datum $u_\phi(0) = \phi$. By Theorem 4.7, after extracting a subsequence, we can decompose the sequence in (6.2) into a sum of Euclidean profiles $\widetilde{\varphi}^{\alpha}_{\mathcal{O}^\alpha,k}$ and a sum of scale-one profiles $\widetilde{\omega}^{\beta}_{\mathcal{O}^\beta,k}$, where $\alpha, k \in \mathbb{N}^+$. More precisely, we can write

$$u_k(0) = \sum_{\alpha=1}^{J} \widetilde{\varphi}^{\alpha}_{\mathcal{O}^\alpha,k} + \sum_{\beta=1}^{J} \widetilde{\omega}^{\beta}_{\mathcal{O}^\beta,k} + R_k^J$$

where $R_k^J$ satisfies

$$\limsup_{J\to\infty} \limsup_{k\to\infty} \|e^{it\Delta} R_k^J\|_{Z(I_k)} = 0. \tag{6.3}$$

Moreover, the almost orthogonality in Proposition 4.7 and the almost orthogonality of nonlinear profiles in Lemma 5.9 give that

$$\lim_{J\to J^*} \Big( \sum_{1\leq\alpha\leq J} \mathcal{L}_E(\alpha) + \sum_{1\leq\beta\leq J} \mathcal{L}_1(\beta) + \lim_{k\to+\infty} L(R_k^J) \Big) = L_{max}$$

where we denote

$$\mathcal{L}_E(\alpha) := \lim_{k\to+\infty} L(\widetilde{\varphi}^{\alpha}_{\mathcal{O}^{\alpha}_k}) \in [0, L_{max}],$$
$$\mathcal{L}_1(\beta) := \lim_{k\to+\infty} L(\widetilde{\omega}^{\beta}_{\mathcal{O}^{\beta}_k}) \in [0, L_{max}]. \tag{6.4}$$

Note that in this section, we will constantly use $\varphi$ and $E$ in the Euclidean profile related context, and employ $\omega$ and 1 in the scale-one profile related context.

Denote $U_k^{\mathrm{E},\alpha}$ to be the solution of (1.1) with initial datum $\widetilde{\varphi}^{\alpha}_{\mathcal{O}^{\alpha},k}$. By Lemma 5.9 and Theorem 4.4, for large enough $k$, $U_k^{\mathrm{E},\alpha}$ uniquely exists also in the time interval $[-T_k, T^k]$ and satisfies that

$$\sup_{t\in[-T_k,T^k]} \|U_k^{\mathrm{E},\alpha}\|_{\dot{H}^1} = \mathcal{L}_{\mathrm{E}}(\alpha) \le \sup_{t\in[-T_k,T^k]} \|u_k\|_{\dot{H}^1} < \|W\|_{\dot{H}^1(\mathbb{R}^4)},$$

which matches the condition (5.11) and hence allows us to apply the decomposition of $U_k^{\mathrm{E},\alpha}$ in Proposition 5.7 later.

**Case 1:** No profiles.
That is, there are neither Euclidean profiles nor scale-one profiles, then

$$u_k(0) = R_k^J.$$

Taking $J$ sufficiently large, using the smallness of $R_k^J$ (6.3), we see that

$$\|e^{it\Delta}u_k(0)\|_{Z([-1,1])} = \|e^{it\Delta}R_k^J\|_{Z([-1,1])} \le \frac{\delta_0}{2}$$

for $k$ sufficiently large, where $\delta_0 = \delta_0(L_{max})$ is given in Theorem 4.3. Then we see from Theorem 4.3 that $u_k$ can be extended on $[-1, 1]$ and that

$$\|u_k\|_{Z([-1,1])} \le \delta_0^{\frac{5}{3}}$$

which contradicts our assumption (6.2).

**Case 2:** Exactly one Euclidean profile.
That is, there is only one Euclidean profile $\widetilde{\varphi}^1_{\mathcal{O}^1,k}$ such that $\mathcal{L}_E(1) = L_{max}$, where $\mathcal{L}_E(1)$ is defined in (6.4). That is,

$$u_k(0) = \widetilde{\varphi}_{\mathcal{O}_k} + o_k(1)$$

in $H^1$, where $\mathcal{O}$ is a Euclidean frame.

In this case, since from Proposition 5.2 the corresponding nonlinear profile $U_k$ satisfies

$$\|U_k\|_{Z([-1,1])} \lesssim_{E_{\mathbb{R}^4}(\varphi)} 1 \text{ and } \lim_{k\to\infty} \|U_k(0) - u_k(0)\|_{H^1} = 0$$

we may use Theorem 4.5 to deduce that

$$\|u_k\|_{Z([-1,1])} \lesssim \|u_k\|_{X^1([-1,1])} \lesssim_{L_{max}} 1$$

with contradicts (6.2).

**Case 3:** Exactly one scale-one profile.
That is, $\mathcal{L}_1(1) = L_{max}$ and

$$u_k(0) = \widetilde{\omega}_{\mathcal{O}_k} + o_k(1) \quad \text{ in } H^1,$$

where $\mathcal{O}$ is a scale-one frame. Hence $u_k(\cdot - x_k, 0) \to \omega$ strongly in $H^1$. We define $U$ is the solution of (1.1) with the initial datum $U(0) = \omega$. First, we consider the case $|T_k| \to 0$ or $|T^k| \to 0$ as $k \to \infty$. For simplicity, we only consider the case $|T^k| \to 0$. By the local well-posedness theory, there exists $\delta > 0$ such that $\|U\|_{Z(0,\delta)} \le 1$. By the stability theorem (Theorem 4.5) and the strong convergence $u_k(\cdot - x_k, 0) \to \omega$, for sufficient large $k$, we have that

$$\|u_k\|_{Z((0,T^k))} \le \|U\|_{Z((0,\delta))} \le 1,$$

which contradicts (6.2).

Second, we consider the case $\limsup_{k\to\infty}|T_k|>0$ and $\limsup_{k\to\infty}|T^k|>0$. Suppose $\mathcal{T}^\infty:=\lim_{k\to\infty}T^k>0$ and $\mathcal{T}_\infty:=\lim_{k\to\infty}T_k>0$ up to a subsequence. Suppose the maximal time of existence of $U$ is $(-T_-,T^+)$. We claim that $(-T_-,T^+)\subset(-\mathcal{T}_\infty,\mathcal{T}^\infty)$. In fact, for arbitrary small $\delta>0$, if $T^+>\mathcal{T}^\infty+2\delta$ by Theorem 4.4 and the strong convergence, for sufficiently large $k$ we have

$$\|u_k\|_{Z((0,T^k))}\lesssim\|U\|_{Z((0,T^+-\delta))}<\infty,$$

which contradicts (6.2). So $T^+\le\mathcal{T}^\infty$ and similarly we have $T_-\le\mathcal{T}_\infty$.

Then using $(-T_-,T^+)\subset(-\mathcal{T}_\infty,\mathcal{T}^\infty)$ and the strong convergence, we have

$$\sup_{t\in(-T_-,T^+)}\|U(t)\|_{\dot H^1}\le\mathcal{L}_1(1)=L_{max}<\|W\|_{\dot H^1(\mathbb{R}^4)}<\infty. \tag{6.5}$$

However, by Theorem 4.4 and the maximal time of existence, $\|U\|_{Z(-T_-,T^+)}=\infty$ which contradicts (6.5).

**Case 4:** At least two of all profiles are nonzero (multiple profiles).

In this case, $\mathcal{L}_\mu(1)<L_{max}$ for $\mu\in\{E,1\}$. We will construct an approximate solution of (1.1) with initial datum $u_k(0)$ whose $Z$-norm is finite and then derive a contradiction using Theorem 4.3.

First, there exists $\eta>0$ such that for all $\alpha\ge1$, $\mu\in\{E,1\}$ and $\mathcal{L}_\mu(\alpha)<L_{max}-\eta$, we have that all nonlinear profiles are global and satisfy, for any $k$, $\alpha\ge1$ and $\mu\in\{E,1\}$ (after extracting a subsequence)

$$\|U_k^{\mu,\alpha}\|_{Z((-2,2))}\le2\Lambda(L_{max}-\eta/2)\lesssim1.$$

Here we note that from now on all implicit constants are allowed to depend on $\Lambda(L_{max}-\eta/2)$. Using Theorem 4.3, it follows that

$$\|U_k^{\mu,\alpha}\|_{X^1((-1,1))}\lesssim1.$$

For $J,k\ge1$, we define

$$U^J_{prof,k}:=\sum_{1\le\alpha\le J}\sum_{\mu\in\{E,1\}}U_k^{\mu,\alpha}.$$

Now we aim to prove that $\|U^J_{prof,k}\|_{X^1((-1,1))}\lesssim1$.

More precisely, we can show that there exists a constant $C\lesssim1$ such that for all $k$ sufficiently large

$$\|U^J_{prof,k}\|^2_{X^1}+\sum_{1\le\alpha\le J}\sum_{\mu\in\{E,1\}}\|U_k^{\mu,\alpha}\|^2_{X^1}\le C^2$$

uniformly in $J$. We take $\delta_0(2L_{max})$ defined in Theorem 4.3. It is obvious that there are finitely many profiles such that $\mathcal{L}_\mu(\alpha)\ge\frac{\delta_0}{2}$, $\mu\in\{E,1\}$. Without loss of generality, we may assume that there exists an $A\in\mathbb{N}^+$ such that for all $\alpha\ge A,\mathcal{L}_\mu(\alpha)\le\delta_0$. Similar to the defocusing case, we have

$$\begin{aligned}\|U^J_{prof,k}\|_{X^1((-1,1))}&=\|\sum_{1\le\alpha\le J}\sum_{\mu\in\{E,1\}}U_k^{\mu,\alpha}\|_{X^1((-1,1))}\\&\le\|\sum_{1\le\alpha\le A}\sum_{\mu\in\{E,1\}}U_k^{\mu,\alpha}\|_{X^1(-1,1))}+\|\sum_{A\le\alpha\le J}\sum_{\mu\in\{E,1\}}(U_k^{\mu,\alpha}-e^{it\Delta}U_k^{\mu,\alpha}(0))\|_{X^1((-1,1))}\\&\quad+\|e^{it\Delta}\sum_{A\le\alpha\le J}\sum_{\mu\in\{E,1\}}U_k^{\mu,\alpha}(0)\|_{X^1((-1,1))}\\&\le\|\sum_{1\le\alpha\le A}\sum_{\mu\in\{E,1\}}U_k^{\mu,\alpha}\|_{X^1((-1,1))}+\sum_{A\le\alpha\le J}\sum_{\mu\in\{E,1\}}\mathcal{L}_\mu(\alpha)^{\frac32}\\&\quad+\|\sum_{A\le\alpha\le J}\sum_{\mu\in\{E,1\}}U_k^{\mu,\alpha}(0)\|_{H^1(\mathbb{R}\times\mathbb{T}^3)}\\&\lesssim A+\sum_{A\le\alpha\le J}\sum_{\mu\in\{E,1\}}\mathcal{L}_\mu(\alpha)^{\frac32}+\|\sum_{A\le\alpha\le J}\sum_{\mu\in\{E,1\}}U_k^{\mu,\alpha}(0)\|_{H^1(\mathbb{R}\times\mathbb{T}^3)}\lesssim1.\end{aligned}$$

The bound on $\sum_{1\le\alpha\le J}\sum_{\mu\in\{E,1\}}\|U_k^{\mu,\alpha}\|_{X^1}^2$ is similar.

Then we are now ready to construct the approximate solution. Let $F(z)=z|z|^2$ and we have

$$F^{'}(G)u=2G\bar{G}u+G^2\bar{u}.$$

For each $B$ and $J$, we define $g_k^{B,J}$ to be the solution of the initial value problem:

$$i\partial_t g+\Delta_{\mathbb{R}\times\mathbb{T}^3}g+F^{'}(U^B_{prof,k})g=0,\quad g(0)=R_k^J.$$

The solution $g_k^{B,J}$ is well defined on $(-1,1)$ for $k>k_0(B,J)$ and satisfies:

$$\|g_k^{B,J}\|_{X^1((-1,1))}\le C^{'},$$

for some $C^{'}$ independent of $J$ and $B$. In fact, this follows from splitting $\mathbb{R}$ into $O(C)$ intervals $I_j$ over which $\|U^B_{prof,k}\|_{Z(I_j)}$ is small and applying the local theory on each subinterval.

Now we can define the approximate solution using the sum of all the nonlinear profiles. We let ($A$ will be chosen shortly)

$$U_k^{app,J}:=U^A_{prof,J}+g_k^{A,J}+U^{>A}_{prof,k}\quad\text{where}\quad U^{>A}_{prof,k}=\sum_{A<\alpha\le J}\sum_{\mu}U_k^{\mu,\alpha}.$$

satisfy the bound for any $1\le A\le J$,

$$\|U_k^{app,J}\|_{X^1((-1,1))}\le 3(C+C^{'})$$

for $k$ large enough (depending on $J$). Applying Theorem 4.3 with $M=6(1+C+C^{'})$ we have a $\varepsilon_1=\varepsilon_1(M)\le\frac{1}{K(1+C+C^{'})}$ for some $K$ sufficiently large, such that if the error term defined below (6.7) with $N$-norm bounded by $\varepsilon_1$. Then we can upgrade the uniform $X^1((-1,1))$ bounds into a uniform bound on $\|u_k\|_{Z((-1,1))}$, hence deriving a contradiction.

First we choose $A$ such that:

$$\|U^{>A}_{prof,k}\|^2_{X^1((-1,1))}+\sum_{A<\alpha\le J}\sum_{\mu}\|U_k^{\mu,\alpha}\|^2_{X^1((-1,1))}\le\varepsilon_1^{10}. \tag{6.6}$$

for any $J\ge A$ and $k$ sufficiently large.

After fixing $A$ we can rewrite the error term as

$$\begin{aligned}
e_k^J&=(i\partial_t+\Delta)U_k^{app,J}-F(U_k^{app,J}) &(6.7)\\
&=-F(U^A_{prof,k}+g_k^{A,J}+U^{>A}_{prof,k})+\sum_{1\le\alpha\le J,\mu}F(U_k^{\mu,\alpha})+F^{'}(U^A_{prof,k})g_k^{A,J}\\
&=-F(U^A_{prof,k}+g_k^{A,J}+U^{>A}_{prof,k})+F(U^A_{prof,k}+g_k^{A,J})+F(U^{>A}_{prof,k}) &(6.8)\\
&\quad-F(U^A_{prof,k}+g_k^{A,J})+F(U^A_{prof,k})+F^{'}(U^A_{prof,k})g_k^{A,J} &(6.9)\\
&\quad-F(U^A_{prof,k})+\sum_{1\le\alpha\le A}F(U_k^{\mu,\alpha}) &(6.10)\\
&\quad-F(U^{>A}_{prof,k})+\sum_{A+1\le\alpha\le J,\mu}F(U_k^{\mu,\alpha}). &(6.11)
\end{aligned}$$

We will estimate the four terms separately.

Before we do the calculation, let us state two useful lemmas.

**Lemma 6.2.** *Assume that $U_k^\alpha,U_k^\beta,U_k^\gamma$ are three nonlinear profiles from the set $\{U_k^{\mu,\alpha}:1\le\alpha\le A,\mu\in\{E,1\}\}$ such that $U_k^\alpha$ and $U_k^\beta$ correspond to orthogonal frames. Then for these nonlinear profiles:*

$$\limsup_{k\to+\infty}\|\widetilde{U}_k^\alpha\widetilde{U}_k^\beta\widetilde{U}_k^\gamma\|_{N((-1,1))}=0$$

*where for* $\delta\in\{\alpha,\beta,\gamma\}$ , $\widetilde{U}_k^\delta\in\{U_k^\delta,\bar{U}_k^\delta\}$.

**Lemma 6.3.** *For any fixed A, it holds that:*

$$\limsup_{J\to\infty}\limsup_{k\to\infty}\|g_k^{A,J}\|_{Z((-1,1))}=0.$$

*Proofs of Lemma 6.2 and Lemma 6.3.* The proofs of Lemma 6.2 and Lemma 6.3 are similar to the proofs of the defocusing analogues, i.e. Lemma 6.3 and Lemma 6.4 in [22], respectively. Note that in particular the decomposition of the nonlinear profiles plays an important role in the proofs of Lemma 6.2 and Lemma 6.3. For the nonlinear Euclidean profiles $U^{E,\alpha}$, we have recovered the decomposition as in Proposition 5.7 and for the nonlinear scale-one profiles $U^{1,\beta}$, the corresponding decomposition can be easily achieved similar as in Corollary 5.8. □

To finish the proof for the last scenario, we discuss the estimations of terms in (6.8), (6.9), (6.10) and (6.11) respectively. The process is almost the same as the defocusing case. Lemma 6.2, Lemma 6.3, and nonlinear estimates in [22] (Lemma 4.2) are often used.

First, we estimate:

$$\|(6.8)\|_{N((-1,1))}\lesssim(\|U^A_{prof,k}+g_k^{A,J}\|_{X^1((-1,1))}+\|U^{>A}_{prof,k}\|_{X^1((-1,1))})^2\|U^{>A}_{prof,k}\|_{X^1((-1,1))}<\frac{\varepsilon_1}{4}$$

for $k$ large enough.

By Lemma 6.3, we estimate:

$$\|(6.9)\|_{N((-1,1))}\lesssim(\|U^A_{prof,k}\|_{X^1((-1,1))}+\|g_k^{A,J}\|_{X^1((-1,1))})^2\|g_k^{A,J}\|_{Z'((-1,1))}\lesssim(Q+Q')^2\|g_k^{A,J}\|_{Z'((-1,1))}<\frac{\varepsilon_1}{4}$$

for $k>k_0(J)$ if $J$ is big enough.

By Lemma 6.2, we estimate:

$$\|(6.10)\|_{N((-1,1))}\lesssim_A\sum_{(\alpha_1,\mu_1)\neq(\alpha_2,\mu_2)}\|\widetilde{U}_k^{\mu_1,\alpha_1}\widetilde{U}_k^{\mu_2,\alpha_2}\widetilde{U}_k^{\mu_3,\alpha_3}\|_{N((-1,1))}<\frac{\varepsilon_1}{4}$$

if $k$ is big enough.

By (6.6), we estimate:

$$\|(6.11)\|_{N((-1,1))}\lesssim\|U^{>A}_{prof,k}\|^3_{X^1((-1,1))}+\sum_{A<\alpha\le J}\|U^{\mu,\alpha}_{prof,k}\|^3_{X^1((-1,1))}\le\frac{\varepsilon_1}{4}.$$

As last, using Theorem 4.3, we obtain that $u_k$ extends as a solution in $X^1((-1,1))$ satisfying:

$$\|u_k\|_{Z((-1,1))}<+\infty$$

which contradicts (6.2).

□

**Acknowledgments.** The authors thank the organizers of *2019 Conference on Nonlinear Partial Differential Equations and Applications* at University of Michigan, where this work began.

## References

[1] T. Aubin. Équations différentielles non linéaires et problème de Yamabe concernant la courbure scalaire. *J. Math. Pures Appl. (9)*, 55(3):269–296, 1976.

[2] A. Barron. On Global-in-Time Strichartz Estimates for the Semiperiodic Schrödinger Equation. *arXiv preprint arXiv:1901.01663*, 2019.

[3] J. Bourgain. Exponential sums and nonlinear Schrödinger equations. *Geom. Funct. Anal.*, 3(2):157–178, 1993.

[4] J. Bourgain. Fourier transform restriction phenomena for certain lattice subsets and applications to nonlinear evolution equations. I. Schrödinger equations. *Geom. Funct. Anal.*, 3(2):107–156, 1993.

[5] J. Bourgain and C. Demeter. The proof of the $l^2$ decoupling conjecture. *Ann. of Math. (2)*, 182(1):351–389, 2015.

[6] T. Cazenave. *Semilinear Schrödinger equations*, volume 10 of *Courant Lecture Notes in Mathematics*. New York University, Courant Institute of Mathematical Sciences, New York; American Mathematical Society, Providence, RI, 2003.

[7] T. Cazenave and F. B. Weissler. The Cauchy problem for the critical nonlinear Schrödinger equation in $H^s$. *Nonlinear Anal.*, 14(10):807–836, 1990.

[8] X. Cheng, Z. Guo, K. Yang, and L. Zhao. On scattering for the cubic defocusing nonlinear schrödinger equation on waveguide $\mathbb{R}^2 \times \mathbb{T}$. *arXiv preprint arXiv:1705.00954*, 2017.

[9] X. Cheng, Z. Guo, and Z. Zhao. On scattering for the defocusing quintic nonlinear schrödinger equation on the two-dimensional cylinder. *arXiv preprint arXiv:1809.01527*, 2018.

[10] X. Cheng, Z. Zhao, and J. Zheng. Well-posedness for energy-critical nonlinear schrödinger equation on waveguide manifold. *arXiv preprint arXiv:1911.00324*, 2019.

[11] J. Colliander, M. Keel, G. Staffilani, H. Takaoka, and T. Tao. Global well-posedness and scattering for the energy-critical nonlinear Schrödinger equation in $\mathbb{R}^3$. *Ann. of Math. (2)*, 167(3):767–865, 2008.

[12] B. Dodson. Global well-posedness and scattering for the defocusing, $L^2$-critical nonlinear Schrödinger equation when $d \geq 3$. *J. Amer. Math. Soc.*, 25(2):429–463, 2012.

[13] B. Dodson. Global well-posedness and scattering for the defocusing, $L^2$-critical, nonlinear Schrödinger equation when $d = 2$. *Duke Math. J.*, 165(18):3435–3516, 2016.

[14] B. Dodson. *Defocusing nonlinear Schrödinger equations*, volume 217 of *Cambridge Tracts in Mathematics*. Cambridge University Press, Cambridge, 2019.

[15] B. Dodson. Global well-posedness and scattering for the focusing, cubic Schrödinger equation in dimension $d = 4$. *Ann. Sci. Éc. Norm. Supér. (4)*, 52(1):139–180, 2019.

[16] Z. Hani and B. Pausader. On scattering for the quintic defocusing nonlinear Schrödinger equation on $\mathbb{R} \times \mathbb{T}^2$. *Comm. Pure Appl. Math.*, 67(9):1466–1542, 2014.

[17] E. Hebey. *Sobolev spaces on Riemannian manifolds*, volume 1635 of *Lecture Notes in Mathematics*. Springer-Verlag, Berlin, 1996.

[18] E. Hebey and M. Vaugon. The best constant problem in the Sobolev embedding theorem for complete Riemannian manifolds. *Duke Math. J.*, 79(1):235–279, 1995.

[19] S. Herr, D. Tataru, and N. Tzvetkov. Global well-posedness of the energy-critical nonlinear Schrödinger equation with small initial data in $H^1(\mathbb{T}^3)$. *Duke Math. J.*, 159(2):329–349, 2011.

[20] S. Herr, D. Tataru, and N. Tzvetkov. Strichartz estimates for partially periodic solutions to Schrödinger equations in $4d$ and applications. *J. Reine Angew. Math.*, 690:65–78, 2014.

[21] A. D. Ionescu and B. Pausader. The energy-critical defocusing NLS on $\mathbb{T}^3$. *Duke Math. J.*, 161(8):1581–1612, 2012.

[22] A. D. Ionescu and B. Pausader. Global well-posedness of the energy-critical defocusing NLS on $\mathbb{R} \times \mathbb{T}^3$. *Comm. Math. Phys.*, 312(3):781–831, 2012.

[23] A. D. Ionescu, B. Pausader, and G. Staffilani. On the global well-posedness of energy-critical Schrödinger equations in curved spaces. *Anal. PDE*, 5(4):705–746, 2012.

[24] C. E. Kenig and F. Merle. Global well-posedness, scattering and blow-up for the energy-critical, focusing, non-linear Schrödinger equation in the radial case. *Invent. Math.*, 166(3):645–675, 2006.

[25] C. E. Kenig and F. Merle. Global well-posedness, scattering and blow-up for the energy-critical focusing non-linear wave equation. *Acta Math.*, 201(2):147–212, 2008.

[26] R. Killip and M. Vişan. Nonlinear Schrödinger equations at critical regularity. In *Evolution equations*, volume 17 of *Clay Math. Proc.*, pages 325–437. Amer. Math. Soc., Providence, RI, 2013.

[27] R. Killip and M. Vişan. Scale invariant Strichartz estimates on tori and applications. *Math. Res. Lett.*, 23(2):445–472, 2016.

[28] R. Killip and M. Visan. The focusing energy-critical nonlinear Schrödinger equation in dimensions five and higher. *Amer. J. Math.*, 132(2):361–424, 2010.
[29] B. Pausader, N. Tzvetkov, and X. Wang. Global regularity for the energy-critical NLS on $\mathbb{S}^3$. *Ann. Inst. H. Poincaré Anal. Non Linéaire*, 31(2):315–338, 2014.
[30] E. Ryckman and M. Vişan. Global well-posedness and scattering for the defocusing energy-critical nonlinear Schrödinger equation in $\mathbb{R}^{1+4}$. *Amer. J. Math.*, 129(1):1–60, 2007.
[31] G. Talenti. Best constant in Sobolev inequality. *Ann. Mat. Pura Appl. (4)*, 110:353–372, 1976.
[32] T. Tao. *Nonlinear dispersive equations*, volume 106 of *CBMS Regional Conference Series in Mathematics*. Published for the Conference Board of the Mathematical Sciences, Washington, DC; by the American Mathematical Society, Providence, RI, 2006. Local and global analysis.
[33] H. Yue. Global Well-Posedness of the Energy-Critical Nonlinear Schrödinger Equation on $\mathbb{T}^4$. *arXiv preprint arXiv:1805.09816*, 2018.
[34] Z. Zhao. On scattering for the defocusing nonlinear schrödinger equation on waveguide $\mathbb{R}^m \times \mathbb{T}$ (when $m = 2, 3$). *arXiv preprint arXiv:1712.01266*, 2017.
[35] Z. Zhao. Global well-posedness and scattering for the defocusing cubic Schrödinger equation on waveguide $\mathbb{R}^2 \times \mathbb{T}^2$. *J. Hyperbolic Differ. Equ.*, 16(1):73–129, 2019.
[36] Z. Zhao and J. Zheng. Long time dynamics for defocusing cubic NLS on three dimensional product space. *arXiv preprint arXiv:2003.00665*, 2020.

Xueying Yu
Department of Mathematics, MIT
77 Massachusetts Ave, Cambridge, MA 02139,

*E-mail address*: xueyingy@mit.edu

Haitian Yue
Department of Mathematics, University of Southern California
3620 S. Vermont Ave, Los Angeles, CA 90089.

*E-mail address*: haitiany@usc.edu

Zehua Zhao
Department of Mathematics, University of Maryland
William E. Kirwan Hall, 4176 Campus Dr. College Park, MD 20742-4015

*E-mail address*: zzh@umd.edu